\documentclass[12pt]{article}

\usepackage[applemac]{inputenc}
\usepackage{amsmath,amssymb,fullpage}
\usepackage{showlabels}
\usepackage{color}

\newcommand{\dproof}{\noindent {Proof.} \quad}
\newcommand{\fproof}{\hfill $\square$ \bigskip}

\newtheorem{definition}{Definition}[section]

\newtheorem{theorem}[definition]{Theorem}
\newtheorem{problem}[definition]{Problem}
\newtheorem{remark}[definition]{ \it Remark}
\newtheorem{coro}[definition]{Corollary}

\numberwithin{equation}{section}

\def\1B{\text{1\!\!I}}

\begin{document}
\date{9 September 2015}
\title{ Stochastic differential games with inside information}

\author{
Olfa Draouil$^{1}$ and Bernt \O ksendal$^{2,3}$}

\footnotetext[1]{Department of Mathematics, University of Tunis El Manar, Tunis, Tunisia.
Email: {\tt olfadraouil@hotmail.fr}}

\footnotetext[2]{Department of Mathematics, University of Oslo, P.O. Box 1053 Blindern, N--0316 Oslo, Norway. Email: {\tt oksendal@math.uio.no}}

\footnotetext[3]{This research was carried out with support of CAS - Centre for Advanced Study, at the Norwegian Academy of Science and Letters, within the research program SEFE.}
\maketitle


\begin{abstract}
We study stochastic differential games of jump diffusions, where the players have access to inside information. Our approach is based on anticipative stochastic calculus, white noise, Hida-Malliavin calculus, forward integrals and the Donsker delta functional. We obtain a characterization of Nash equilibria of such games in terms of the corresponding Hamiltonians. This is used to study applications to insider games in finance, specifically optimal insider consumption and optimal insider portfolio under model uncertainty.
\end{abstract}

\section{Introduction}
In this paper we present a general method for solving \emph{optimal insider games}, i.e. optimal insider games problems where the players has access to some future information about the system. This inside information in the control processes puts the problem outside the context of semimartingale theory, and we therefore apply general \emph{anticipating white noise calculus}, including \emph{forward integrals} and \emph{Hida-Malliavin calculus}. Combining this with the \emph{Donsker delta functional} for the random variable $Y=(Y_1,Y_2)$ which represents the inside information, we are able to prove both a sufficient and a necessary maximum principle for the optimal insider games of such systems.

We now explain this in more detail:\\

The system we consider, is described by a stochastic differential equation driven by a Brownian motion $B(t)$ and an independent compensated Poisson random measure $\tilde{N}(dt,d\zeta)$, jointly defined on a filtered probability space $(\Omega, \mathbb{F}=\{ \mathcal{F}_t \}_{t \geq 0},\mathbf{P})$ satisfying the usual conditions. We assume that the inside information is of \emph{initial enlargement} type. Specifically, we assume  that the two inside filtrations $\mathbb{H}^1,\mathbb{H}^2$ representing the information flows available to player 1 and player 2, respectively,  have the form

\begin{equation}\label{eq1.1}
 \mathbb{H}^i= \{ \mathcal{H}_t^i\}_{t \geq 0}, \text{ where } \mathcal{H}^i_t = \mathcal{F}_t \vee Y_i, \quad i=1,2
\end{equation}
for all $t$, where $Y_i$ is a given $\mathcal{F}_{T_0}$-measurable random variable, for some fixed $T_0 > T >t$.
Here the insider control process $u(t)=(u_1(t),u_2(t))$, where $u_i(t)$ is the control of player i; i=1,2.
Thus we assume that the value at time $t$ of our insider control process $u_i(t)$ is allowed to depend on both $Y_i$ and $\mathcal{F}_t; i=1,2.$ In other words, $u_i$ is assumed to be $\mathbb{H}^i$-adapted for $i=1,2$. Therefore they  have the form
\begin{equation}\label{eq1.2}
    u_i(t,\omega) = \overline{u}_i(t, Y_i, \omega)
\end{equation}
for some function $\overline{u}_i : [0, T]\times \mathbb{R}\times\Omega\rightarrow \mathbb{R}$
such that $\overline{u}_i(t, y_i)$ is $\mathbb{F}$-adapted for each $y_i \in\mathbb{R}$. For simplicity (albeit with some abuse of notation) we will in the following write $u_i$ in stead of $\overline{u}_i; i=1,2$.
Consider a controlled stochastic process $X(t)=X^u(t)$ of the form
\begin{equation}\label{eq2.1}
    \left\{
  \begin{array}{l}
dX(t)=dX^u(t)=b(t,X(t),u_1(t),u_2(t),Y_1,Y_2)dt+\sigma(t,X(t),u_1(t),u_2(t),Y_1,Y_2)dB(t)\\
+\int_{\mathbb{R}} \gamma(t,X(t),u_1(t),u_2(t),Y_1,Y_2,\zeta)\tilde{N}(dt,d\zeta); \quad t\geq 0\\
X(0)=x(Y)\in \mathbb{R},
  \end{array}
    \right.
\end{equation}
where
$u_i(t)=u_i(t,y_i)_{y_i=Y_i}$  is the control process of insider $i; i=1,2,$ and the (anticipating) stochastic integrals are interpreted as \emph{forward integrals}, as introduced in \cite{RV} (Brownian motion case) and in \cite{DMOP1} (Poisson random measure case). A motivation for using forward integrals in the modelling of insider control is given in \cite{BO}.
Let $\mathcal{A}_i$  denote a given set of admissible $\mathbb{H}^i-$adapted controls $u_i$ of player $i$, with values in $\mathbf{A}_i\subset\mathbb{R}^d, d\geq1; i=1,2$. Denote $\mathbb{U}=\mathbf{A}_1\times\mathbf{A}_2$.
 Then $X(t)$ is $\mathbb{F}\vee Y_1 \vee Y_2$-adapted.
The \emph{performance functional} $J_i(u); u=(u_1,u_2)$ of player $i$ is defined by, writing $y=(y_1,y_2)$ and $dy=dy_1dy_2$,
\begin{eqnarray}\label{performance}
  J_i(u) &=& \mathbb{E}[\int_0^T f_i(t, X(t),u_1(t),u_2(t),Y)dt +g_i(X(T),Y)]   \nonumber\\
   &=& \mathbb{E}[\int_{\mathbb{R}^2}\{\int_0^T f_i(t, x(t,y_1,y_2),u_1(t,y_1),u_2(t,y_2),y)\mathbb{E}[\delta_{Y_1,Y_2}(y_1,y_2)|\mathcal{F}_t]dt \nonumber\\
   &+& g_i(x(T,y_1,y_2),y_1,y_2)\mathbb{E}[\delta_{Y_1,Y_2}(y_1,y_2)|\mathcal{F}_T]\}dy] ; \quad i=1,2.
\end{eqnarray}
A \emph{Nash equilibrium} for the game (\ref{eq2.1})-(\ref{performance}) is a pair $\hat{u}=(\hat{u}_1,\hat{u}_2)\in \mathcal{A}_1\times\mathcal{A}_2$ such that
\begin{equation}\label{eq2.10}
\sup_{u_1\in\mathcal{A}_1}J_1(u_1,\hat{u}_2)\leq J_1(\hat{u}_1,\hat{u}_2)
\end{equation}
and
\begin{equation}\label{eq2.11}
\sup_{u_2\in\mathcal{A}_2}J_2(\hat{u}_1,u_2)\leq J_2(\hat{u}_1,\hat{u}_2).
\end{equation}
We use the Donsker delta functional of $Y=(Y_1,Y_2)$ to find a Nash equilibrium for the game (\ref{eq2.1})-(\ref{performance}).

Here is an outline of the content of the paper:\\
\begin{itemize}
\item
In Section $2$ we define the Donsker delta functional.
\item
In Section $3$ we present the general insider optimal control problem for the
stochastic differential games.
\item
In Sections $4$ and $5$  we present a sufficient and a necessary maximum principle, respectively, for the insider game problem.
\item
In Section $6$  we present the zero-sum game case where we distinguish two situations: Situation 1: Both players are still maximizing their own performance functional and situation 2: One of the players is maximizing and the other
is minimizing the performance functional and we write the sufficient and necessary maximum principle corresponding for each situation.
\item
Then in Section $7$ we illustrate our results by applying them to optimal insider consumption and optimal insider portfolio under model uncertainty.
\end{itemize}

\section{The Donsker delta functional}
\begin{definition}
Let $(Y_1,Y_2) :\Omega^2\rightarrow\mathbb{R}$ be a pair of random variables which also belongs to $(\mathcal{S})^{\ast^2}$. Then a continuous functional
\begin{equation}\label{donsker}
    \delta_{Y_1,Y_2}(.): \mathbb{R}\times\mathbb{R}\rightarrow (\mathcal{S})^{\ast}
\end{equation}
is called a Donsker delta functional of $(Y_1,Y_2)$ if it has the property that
\begin{equation}\label{donsker property }
    \int_{\mathbb{R}^2}g(y_1,y_2)\delta_{Y_1,Y_2}(y_1,y_2)dy_1dy_2 = g(Y_1,Y_2) \quad a.s.
\end{equation}
for all (measurable) $g : \mathbb{R}^2 \rightarrow \mathbb{R}$ such that the integral converges.
\end{definition}
For more information about the Donsker delta function and some explicit formulas for it, see \cite{DrO}.

\section{The general insider optimal control problem for the stochastic differential games}
In this section, we formulate and prove a sufficient and a necessary maximum principle for general  stochastic differential games (not necessarily zero-sum games) for insiders. The system we consider, is described by a stochastic differential equation driven by a Brownian motion $B(t)$ and an independent compensated Poisson random measure $\tilde{N}(dt,d\zeta)$, jointly defined on a filtered probability space $(\Omega, \mathbb{F}=\{ \mathcal{F}_t \}_{t \geq 0},\mathbf{P})$ satisfying the usual conditions. We assume that the inside information is of \emph{initial enlargement} type. Specifically, we assume  that the two inside filtrations $\mathbb{H}^1,\mathbb{H}^2$ representing the information flows available to player 1 and player 2, respectively,  have the form

\begin{equation}\label{eq1.1}
 \mathbb{H}^i= \{ \mathcal{H}_t^i\}_{t \geq 0}, \text{ where } \mathcal{H}^i_t = \mathcal{F}_t \vee Y_i, \quad i=1,2
\end{equation}
for all $t$, where $Y_i$ is a given $\mathcal{F}_{T}$-measurable random variable, for some fixed $T>t$.
Here the insider control process $u(t)=(u_1(t),u_2(t))$, where $u_i(t)$ is the control of player i; i=1,2.
Thus we assume that the value at time $t$ of our insider control process $u_i(t)$ is allowed to depend on both $Y_i$ and $\mathcal{F}_t; i=1,2.$ In other words, $u_i$ is assumed to be $\mathbb{H}^i$-adapted for $i=1,2$. Therefore they  have the form
\begin{equation}\label{eq1.2}
    u_i(t,\omega) = \overline{u}_i(t, Y_i, \omega)
\end{equation}
for some function $\overline{u}_i : [0, T]\times \mathbb{R}\times\Omega\rightarrow \mathbb{R}$
such that $\overline{u}_i(t, y_i)$ is $\mathbb{F}$-adapted for each $y_i \in\mathbb{R}$. For simplicity (albeit with some abuse of notation) we will in the following write $u_i$ in stead of $\overline{u}_i; i=1,2$.

Consider a controlled stochastic process $X(t)=X^u(t)$ of the form
\begin{equation}\label{eq2.1}
    \left\{
  \begin{array}{l}
dX(t)=dX^u(t)=b(t,X(t),u_1(t),u_2(t),Y_1,Y_2)dt+\sigma(t,X(t),u_1(t),u_2(t),Y_1,Y_2)dB(t)\\
+\int_{\mathbb{R}} \gamma(t,X(t),u_1(t),u_2(t),Y_1,Y_2,\zeta)\tilde{N}(dt,d\zeta); \quad t\geq 0\\
X(0)=x,\quad\hbox{$x\in \mathbb{R}$,}
  \end{array}
    \right.
\end{equation}
where
$u_i(t)=u_i(t,y_i)_{y_i=Y_i}$  is the control process of insider $i; i=1,2,$ and the (anticipating) stochastic integrals are interpreted as \emph{forward integrals}, as introduced in \cite{RV} (Brownian motion case) and in \cite{DMOP1} (Poisson random measure case). A motivation for using forward integrals in the modelling of insider control is given in \cite{BO}.
Let $\mathcal{A}_i$  denote a given set of admissible $\mathbb{H}^i-$adapted controls $u_i$ of player $i$, with values in $\mathbf{A}_i\subset\mathbb{R}^d, d\geq1; i=1,2$. Denote $\mathbb{U}=\mathbf{A}_1\times\mathbf{A}_2$.
 Then $X(t)$ is $\mathbb{F}\vee Y_1 \vee Y_2$-adapted, and hence using the definition of the Donsker delta functional $\delta_{Y_1,Y_2}(y_1,y_2)$ of $(Y_1,Y_2)$ we get
\begin{equation}\label{eq6}
X(t)=x(t,Y_1,Y_2)=x(t,y_1,y_2)_{y_1=Y_1,y_2=Y_2}=\int_{\mathbb{R}^2}x(t,y_1,y_2)\delta_{Y_1,Y_2}(y_1,y_2)dy_1dy_2
\end{equation}
for some $y_1,y_2$-parametrized process $x(t,y_1,y_2)$ which is $\mathbb{F}$-adapted for each $y_1,y_2$.
Then, again by the definition of the Donsker delta functional and the properties of forward integration (\cite{DrO}), we can write
\begin{align}\label{eq7}
X(t)&= x +\int_0^t b(s,X(s),u_1(s),u_2(s),Y_1,Y_2)ds + \int_0^t \sigma(s,X(s),u_1(s),u_2(s),Y_1,Y_2)dB(s)\nonumber\\
&+\int_0^t \int_{\mathbb{R}} \gamma(s,X(s),u_1(s),u_2(s),Y_1,Y_2,\zeta)\tilde{N}(ds,d\zeta)\nonumber\\
&= x+\int_0^t b(s,x(s,Y_1,Y_2),u_1(s,Y_1),u_2(s,Y_2),Y_1,Y_2)ds \nonumber\\
&+ \int_0^t \sigma(s,x(s,Y_1,Y_2),u_1(s,Y_1),u_2(s,Y_2),Y_1,Y_2)dB(s)\nonumber\\
&+\int_0^t \int_{\mathbb{R}} \gamma(s,x(s,Y_1,Y_2),u_1(s,Y_1),u_2(s,Y_2),Y_1,Y_2,\zeta)\tilde{N}(ds,d\zeta)\nonumber\\
&= x+\int_0^t b(s,x(s,y_1,y_2),u_1(s,y_1),u_2(s,y_2),y_1,y_2)_{y_1=Y_1,y_2=Y_2}ds \nonumber\\
&+ \int_0^t \sigma(s,x(s,y_1,y_2),u_1(s,y_1),u_2(s,y_2),y_1,y_2)_{y_1=Y_1,y_2=Y_2}dB(s)\nonumber\\
&+\int_0^t \int_{\mathbb{R}} \gamma(s,x(s,y_1,y_2),u_1(s,y_1),u_2(s,y_2),y_1,y_2,\zeta)_{y_1=Y_1,y_2=Y_2}\tilde{N}(ds,d\zeta)\nonumber\\
&= x+\int_0^t \int_{\mathbb{R}^2}b(s,x(s,y_1,y_2),u_1(s,y_1),u_2(s,y_2),y_1,y_2)\delta_{Y_1,Y_2}(y_1,y_2)dy_1dy_2ds \nonumber\\
&+ \int_0^t \int_{\mathbb{R}^2}\sigma(s,x(s,y_1,y_2),u_1(s,y_1),u_2(s,y_2),y_1,y_2)\delta_{Y_1,Y_2}(y_1,y_2)dy_1dy_2dB(s)\nonumber\\
&+\int_0^t \int_{\mathbb{R}^2}\int_{\mathbb{R}} \gamma(s,x(s,y_1,y_2),u_1(s,y_1),u_2(s,y_2),y_1,y_2,\zeta)\delta_{Y_1,Y_2}(y_1,y_2)dy_1dy_2\tilde{N}(ds,d\zeta)\nonumber\\
&= x+\int_{\mathbb{R}^2} [\int_0^t b(s,x(s,y_1,y_2),u_1(s,y_1),u_2(s,y_2),y_1,y_2)ds\nonumber\\
&+\int_0^t \sigma(s,x(s,y_1,y_2),u_1(s,y_1),u_2(s,y_2),y_1,y_2)dB(s)\nonumber\\
&+\int_0^t \int_{\mathbb{R}} \gamma(s,x(s,y_1,y_2),u_1(s,y_1),u_2(s,y_2),y_1,y_2,\zeta)\tilde{N}(ds,d\zeta)]\delta_{Y_1,Y_2}(y_1,y_2)dy_1dy_2
\end{align}
Comparing \eqref{eq6} and \eqref{eq7} we see that  \eqref{eq6} holds if we choose $x(t,y)$ for each $y=(y_1,y_2) $ as the solution of the classical SDE
\begin{align}\label{eq8}
&dx(t,y_1,y_2) = b(t,x(t,y_1,y_2),u_1(t,y_1),u_2(t,y_2),y_1,y_2)dt \nonumber\\
&+ \sigma(t,x(t,y_1,y_2),u_1(t,y_1),u_2(t,y_2),y_1,y_2)dB(t)\nonumber\\
&+ \int_{\mathbb{R}} \gamma(t,x(t,y_1,y_2),u_1(t,y_1),u_2(t,y_2),y_1,y_2,\zeta)\tilde{N}(dt,d\zeta); \quad t\geq 0\\
 & x(0,y)  = x,\quad\hbox{$x\in \mathbb{R}$,}
\end{align}
The \emph{performance functional} $J_i(u); u=(u_1,u_2)$ of player $i$ is defined by
\begin{eqnarray}\label{performance}
  J_i(u) &=& \mathbb{E}[\int_0^T f_i(t, X(t),u_1(t),u_2(t))dt +g_i(X(T))]   \nonumber\\
   &=& \mathbb{E}[\int_{\mathbb{R}^2}\{\int_0^T f_i(t, x(t,y_1,y_2),u_1(t,y_1),,u_2(t,y_2),y_1,y_2)\mathbb{E}[\delta_{Y_1,Y_2}(y_1,y_2)|\mathcal{F}_t]dt \nonumber\\
   &+& g_i(x(T,y_1,y_2),y_1,y_2)\mathbb{E}[\delta_{Y_1,Y_2}(y_1,y_2)|\mathcal{F}_T]\}dy] ; \quad i=1,2.
\end{eqnarray}
A \emph{Nash equilibrium} for the game (\ref{eq2.1})-(\ref{performance}) is a pair $\hat{u}=(\hat{u}_1,\hat{u}_2)\in \mathcal{A}_1\times\mathcal{A}_2$ such that
\begin{equation}\label{eq2.10}
\sup_{u_1\in\mathcal{A}_1}J_1(u_1,\hat{u}_2)\leq J_1(\hat{u}_1,\hat{u}_2)
\end{equation}
and
\begin{equation}\label{eq2.11}
\sup_{u_2\in\mathcal{A}_2}J_2(\hat{u}_1,u_2)\leq J_2(\hat{u}_1,\hat{u}_2).
\end{equation}
\section{A sufficient maximum principle}
The problem (\ref{eq2.10})-(\ref{eq2.11}) is a stochastic differential game with a standard (albeit parametrized) stochastic differential equation \eqref{eq8} for the state process $x(t,y_1,y_2)$, but with a non-standard performance functional given by (\ref{performance}). We can solve this problem by a modified maximum principle approach, as follows:\\

Define the \emph{Hamiltonians}
 $ H_i:[0,T]\times\mathbb{R}\times\mathbb{R}\times\mathbb{R}\times\mathbb{U}\times\mathbb{R}\times\mathbb{R}\times\mathcal{R} \times \Omega \rightarrow \mathbb{R}$ by
\begin{align}\label{eq11}
&H_i(t,x,y_1,y_2,u_1,u_2,p,q,r)=H(t,x,y_1,y_2,u_1,u_2,p,q,r,\omega)\nonumber\\
&=\mathbb{E}[\delta_{Y_1,Y_2}(y_1,y_2)|\mathcal{F}_t] f_i(t,x,u_1,u_2,y_1,y_2)+b(t,x,u_1,u_2,y_1,y_2)p \nonumber\\
&+ \sigma(t,x,u_1,u_2,y_1,y_2)q+\int_{\mathbb{R}}\gamma(t,x,u_1,u_2,y_1,y_2)r(t,\zeta)\nu(d\zeta); i=1,2.
\end{align}
Here $\mathcal{R}$ denotes the set of all functions $r(.) : \mathbb{R}\rightarrow  \mathbb{R}$
such that the last integral above converges.
For $i=1,2$ we define the \emph{adjoint} processes $p_i(t,y_1,y_2),q_i(t,y_1,y_2), r_i(t,y_1,y_2,\zeta)$ as the solution of the $y_1,y_2$-parametrised BSDEs
\begin{equation}\label{eq12}
    \left\{
\begin{array}{l}
    dp_i(t,y_1,y_2) = -\frac{\partial H_i}{\partial x}(t,y_1,y_2)dt + q_i(t,y_1,y_2)dB(t)+\int_{\mathbb{R}}r_i(t,y_1,y_2,\zeta)\tilde{N}(dt,d\zeta); \quad 0 \leq t\leq T\\
    p_i(T,y_1,y_2)  =  g'_i(x(T,y_1,y_2),y_1,y_2) \mathbb{E}[\delta_{Y_1,Y_2}(y_1,y_2)|\mathcal{F}_T]
\end{array}
    \right.
\end{equation}

Let $J_i(u(.,y_1,y_2))$ be defined by
\begin{align}\label{J(u)2}
    J_i(u(.,y_1,y_2))&=\mathbb{E}[\int_0^T f_i(t, x(t,y_1,y_2),u_1(t,y_1),u_2(t,y_2),y_1,y_2)\mathbb{E}[\delta_{Y_1,Y_2}(y_1,y_2)|\mathcal{F}_t]dt \nonumber\\
    &+g_i(x(T,y_1,y_2),y_1,y_2)\mathbb{E}[\delta_{Y_1,Y_2}(y_1,y_2)|\mathcal{F}_T]]
\end{align}
Then we see that
\begin{align}
J_i(u_1,u_2)= \int_{\mathbb{R}} \int_{\mathbb{R}} J_i(u(.,y_1,y_2) )dy_1 dy_2
\end{align}

The insider game problem can therefore be written as follows:

\begin{problem}
 For each $y_1,y_2\in \mathbb{R}$, find $(u_1^{\star}(.,y_1),u_2^{\star}(.,y_2)) \in\mathcal{A}_1\times\mathcal{A}_2$ such that
 \begin{equation}\label{eq2.15}
\sup_{u_1(.,y_1)\in\mathcal{A}_1} \int_{\mathbb{R}} \int_{\mathbb{R}} J_1(u_1(.,y_1),u^*_2(.,y_2))dy_1 dy_2\leq \int_{\mathbb{R}} \int_{\mathbb{R}}J_1(u^*_1(.,y_1),u^*_2(.,y_2))dy_1 dy_2
\end{equation}
and
\begin{equation}\label{eq2.16}
\sup_{u_2(.,y_2)\in\mathcal{A}_2}\int_{\mathbb{R}} \int_{\mathbb{R}}J_2(u^*_1(.,y_1),u_2(.,y_2))dy_1 dy_2\leq \int_{\mathbb{R}} \int_{\mathbb{R}}J_2(u^*_1(.,y_1),u*_2(.,y_2)) dy_1 dy_2
\end{equation}
\end{problem}
To study this problem we present two maximum principles for the corresponding games. The first is the following:
\begin{theorem}{[Sufficient maximum principle]}\label{sufficient theorem}\\
Let $(\hat{u_1},\hat{u_2}) \in \mathcal{A}_1\times\mathcal{A}_2$ with associated solution $\hat{x}(t,y_1,y_2),\hat{p}_i(t,y_1,y_2),\hat{q}_i(t,y_1,y_2),\hat{r}_i(t,y_1,y_2,\zeta)$ of \eqref{eq8} and \eqref{eq12}; i=1,2. Assume that the following hold:
\begin{enumerate}
 \item $x \rightarrow g_i(x)$ is concave; $i=1,2$
 \item
 The functions
 \begin{equation}
 \hat{\mathcal{H}}_1(x)= \sup_{u_1\in\mathcal{A}_1}\int_{\mathbb{R}}H_1(t,x,y_1,y_2,u_1,\hat{u_2}(t,y_2),\widehat{p}_1(t,y_1,y_2),\widehat{q}_1(t,y_1,y_2),\hat{r}_1(t,y_1,y_2,\cdot))dy_2
 \end{equation}
 and
 \begin{equation}
 \hat{\mathcal{H}}_2(x)= \sup_{u_2\in\mathcal{A}_2}\int_{\mathbb{R}}H_2(t,x,y_1,y_2,\hat{u}_1(t,y_1),u_2,\widehat{p}_2(t,y_1,y_2),\widehat{q}_2(t,y_1,y_2),\hat{r}_2(t,y_1,y_2,\cdot))dy_1
 \end{equation}
 are concave for all $t,y_1,y_2$

\item
 \begin{align}
 \sup_{u_1\in\mathbf{A}_1}&\int_{\mathbb{R}}H_1\big(t,\widehat{x}(t,y_1,y_2),u_1,\hat{u_2}(t,y_2),\widehat{p}_1(t,y_1,y_2),\widehat{q}_1(t,y_1,y_2),\hat{r}_1(t,y_1,y_2,\cdot)\big)dy_2\nonumber\\
      &=\int_{\mathbb{R}}H_1\big(t,\widehat{x}(t,y_1,y_2),\widehat{u}_1(t,y_1),\hat{u_2}(t,y_2),\widehat{p}_1(t,y_1,y_2),\widehat{q}_1(t,y_1,y_2),\hat{r}_1(t,y_1,y_2,\cdot)\big)dy_2\nonumber\\
      & \text{ for all } t,y_1.
\end{align}
\item
 \begin{align}
 \sup_{u_2\in\mathbf{A}_2}&\int_{\mathbb{R}}H_2\big(t,\widehat{x}(t,y_1,y_2),\hat{u}_1(t,y_1),u_2,\widehat{p}_2(t,y_1,y_2),\widehat{q}_2(t,y_1,y_2),\hat{r}_2(t,y_1,y_2,\cdot)\big)dy_1\nonumber\\
  &=\int_{\mathbb{R}}H_2\big(t,\widehat{x}(t,y_1,y_2),\widehat{u}_1(t,y_1),\hat{u_2}(t,y_2),\widehat{p}_2(t,y_1,y_2),\widehat{q}_2(t,y_1,y_2),\hat{r}_2(t,y_1,y_2,\cdot)\big) dy_1\nonumber\\
  &\text{ for all } t,y_2.
\end{align}

 \end{enumerate}
Then $(u^*_1(.,y_1),u^*_2(.,y_2)) := (\widehat{u}_1(.,y_1),\widehat{u}_2(.,y_2))$ is a Nash equilibrium for the problem \eqref{eq2.15}-\eqref{eq2.16}.
\end{theorem}
\dproof  By considering an increasing sequence of stopping times $\tau_n$ converging to $T$, we may assume that all local integrals appearing in the computations below are martingales and have expectation 0.  We omit the details in this argument. See \cite{OS2}.\\
We first prove that
\begin{equation}\label{eq2.17}
\sup_{u_1(.,y_1)\in\mathcal{A}_1}\int_{\mathbb{R}}\int_{\mathbb{R}}J_1(u_1(.,y_1),\hat{u}_2(.,y_2))dy_1dy_2\leq \int_{\mathbb{R}}\int_{\mathbb{R}}J_1(\hat{u}_1(.,y_1),\hat{u}_2(.,y_2))dy_1dy_2
\end{equation}
Choose arbitrary $u_1(.,y_1)\in\mathcal{A}_1$ and let us in the following, for simplicity of notation, put  $x(t,y_1,y_2)=x^{u_1,\hat{u}_2}(t,y_1,y_2), \hat{x}(t,y_1,y_2)=x^{\hat{u}_1,\hat{u}_2}(t,y_1,y_2),\\
 b(t,y_1,y_2)=b(t,x(t,y_1,y_2),u_1(t,y_1),u_2(t,y_2),\omega),\hat{b}(t,y_1,y_2)=b(t,\hat{x}(t,y_1,y_2),u_1(t,y_1),\hat{u}_2(t,y_2),\omega)$
and similarly with $\sigma(t,y_1,y_2), \hat{\sigma}(t,y_1,y_2),\gamma(t,y_1,y_2,\zeta)$, $\hat{\gamma}(t,y_1,y_2,\zeta)$ and $\tilde{x}(t,y_1,y_2)=x(t,y_1,y_2)-\hat{x}(t,y_1,y_2)$.
Let us also put
\begin{equation}
H_1(t,y_1,y_2)=H_1(t,x(t,y_1,y_2),y_1,y_2,u_1(t,y_1),\hat{u_2}(t,y_2),\widehat{p}_1(t,y_1,y_2),\widehat{q}_1(t,y_1,y_2),\hat{r}_1(t,y_1,y_2,\cdot))
\end{equation}
and
\begin{equation}
\hat{H}_1(t,y_1,y_2)=H_1(t,\hat{x}(t,y_1,y_2),y_1,y_2,\hat{u}_1(t,y_1),\hat{u_2}(t,y_2),\widehat{p}_1(t,y_1,y_2),\widehat{q}_1(t,y_1,y_2),\hat{r}_1(t,y_1,y_2,\cdot))
\end{equation}

 Consider
 \begin{equation*}
    \int_{\mathbb{R}}\int_{\mathbb{R}}[J_1(u_1(.,y_1),\hat{u_2}(.,y_2))-J_1(\widehat{u}_1(.,y_1),\hat{u_2}(.,y_2))]dy_1dy_2=I_1+ I_2,
 \end{equation*}
 where
\begin{align}\label{I_1}
I_1= \int_{\mathbb{R}}\int_{\mathbb{R}}\mathbb{E}\big[\int_0^T\{f_1(t,x(t,y_1,y_2),u_1(t,y_1),\hat{u}_2(t,y_2))-f_1(t,\hat{x}(t,y_1,y_2),\hat{u}_1(t,y_1),\hat{u}_2(t,y_2))\}\nonumber\\
\mathbb{E}[\delta_{Y_1,Y_2}(y_1,y_2)|\mathcal{F}_t]dt\big]dy_1dy_2
\end{align}
    and
 \begin{equation}\label{I_2}
     I_2= \int_{\mathbb{R}}\int_{\mathbb{R}}\mathbb{E}\big[\{g_1(x(T,y_1,y_2))-g_1(\widehat{x}(T,y_1,y_2))\}\mathbb{E}[\delta_{Y_1,Y_2}(y_1,y_2)|\mathcal{F}_T]\big]dy_1dy_2.
 \end{equation}
  By the definition of $H_1$ we have
  \begin{eqnarray}\label{II1}
    I_1 &=&  \int_{\mathbb{R}}\int_{\mathbb{R}}\mathbb{E}[\int_0^T\{H_1(t,y_1,y_2)-\widehat{H}_1(t,y_1,y_2)-\widehat{p}_1(t,y_1,y_2)\widetilde{b}_1(t,y_1,y_2) - \widehat{q}_1(t,y_1,y_2)\widetilde{\sigma}(t,u_1,u_2,y_1,y_2)\nonumber\\
    &-&\int_{\mathbb{R}}\hat{r}_1(t,y_1,y_2,\zeta)\tilde{\gamma}(t,y_1,y_2,\zeta)\nu(d\zeta)\}dt]dy_1dy_2.
  \end{eqnarray}
Since $g_1$ is concave we have
 \begin{align}
 \label{II_2}
   I_2 &\leq \int_{\mathbb{R}}\int_{\mathbb{R}}\mathbb{E}[g_1'(\widehat{x}(T,y_1,y_2))\mathbb{E}[\delta_{Y_1,Y_2}(y_1,y_2)|\mathcal{F}_T]\widetilde{x}(T,y_1,y_2)]dy_1dy_2\nonumber\\
   &= \int_{\mathbb{R}}\int_{\mathbb{R}}\mathbb{E}[\widehat{p}_1(T,y_1,y_2)\widetilde{x}(T,y_1,y_2)]dy_1dy_2 \nonumber  \\
    &= \int_{\mathbb{R}}\int_{\mathbb{R}}\mathbb{E}[\int_0^T \widehat{p}_1(t,y_1,y_2) d\widetilde{x}(t,y_1,y_2)+\int_0^T\widetilde{x}(t,y_1,y_2)d\widehat{p}_1(t,y_1,y_2)+\int_0^Td[\hat{p}_1 , \tilde{x}]_t]dy_1dy_2  \nonumber\\
    &= \int_{\mathbb{R}}\int_{\mathbb{R}}\mathbb{E}\big[\int_0^T \widehat{p}_1(t,y_1,y_2) (\widetilde{b}(t,y_1,y_2)dt+\widetilde{\sigma}(t,y_1,y_2)dB(t)+\int_{\mathbb{R}}\tilde{\gamma}(t,y_1,y_2,\zeta)\tilde{N}(dt,d\zeta)) \nonumber\\
    &- \int_0^T\frac{\partial \widehat{H}_1}{\partial x}(t,y_1,y_2)\widetilde{x}(t,y_1,y_2)dt+\int_0^T\widehat{q}_1(t,y_1,y_2)\widetilde{x}(t,y_1,y_2)dB(t)\\ \nonumber
    &+\int_0^T\int_{\mathbb{R}}\widetilde{x}(t,y_1,y_2)\hat{r}_1(t,y_1,y_2,\zeta)\tilde{N}(dt,d\zeta)+\int_0^T\widetilde{\sigma}(t,y_1,y_2)\widehat{q}_1(t,y_1,y_2)dt\\ \nonumber
    &+\int_0^T\int_{\mathbb{R}}\tilde{\gamma}(t,y_1,y_2,\zeta)\hat{r}_1(t,y_1,y_2,\zeta)\nu(d\zeta)dt+\int_0^T\int_{\mathbb{R}}\tilde{\gamma}(t,y_1,y_2,\zeta)\hat{r}_1(t,y_1,y_2,\zeta)\tilde{N}(dt,d\zeta)\big]dy_1dy_2\\  \nonumber
      &= \int_{\mathbb{R}}\int_{\mathbb{R}}\mathbb{E}\big[\int_0^T \widehat{p}_1(t,y_1,y_2) \widetilde{b}(t,y_1,y_2)dt -\int_0^T\frac{\partial \widehat{H}_1}{\partial x}(t,y_1,y_2)\widetilde{x}(t,y_1,y_2)dt
    +\int_0^T\widetilde{\sigma}(t,y_1,y_2)\widehat{q}_1(t)dt \nonumber\\
    &+\int_0^T\int_{\mathbb{R}}\tilde{\gamma}(t,y_1,y_2,\zeta)\hat{r}_1(t,y_1,y_2,\zeta)\nu(d\zeta)dt\big]dy_1dy_2.\nonumber
     \end{align}
Adding (\ref{II1}) - (\ref{II_2}) we get, by concavity of $H_1$,
\begin{align}\label{eq2.28}
  & \int_{\mathbb{R}}\int_{\mathbb{R}}[J_1(u_1(.,y_1),\hat{u}_2(.,y_2))-J_1(\widehat{u}_1(.,y_1),\hat{u}_2(.,y_2))]dy_1dy_2\nonumber\\
  &\leq \int_{\mathbb{R}}\int_{\mathbb{R}}\mathbb{E}\big[\int_0^T \{H_1(t,y_1,y_2)-\widehat{H}_1(t,y_1,y_2)-\frac{\partial \hat{H}_1}{\partial x}(t,y_1,y_2)\widetilde{x}(t,y_1,y_2)\}dt\big]dy_1dy_2
\end{align}
Since $\hat{\mathcal{H}}_1(x)$ is concave, it follows by a standard separating hyperplane argument that there exists a supergradient $a\in\mathbb{R}$ for
$\hat{\mathcal{H}}_1(x)$ at $x=\hat{x}(t,y_1,y_2)$ such that if we define
\begin{equation}
\phi(x)=\hat{\mathcal{H}}_1(x)-\hat{\mathcal{H}}_1(\hat{x}(t,y_1,y_2))-a(x-\hat{x}(t,y_1,y_2))
\end{equation}
then
\begin{equation}
\phi(x)\leq 0 \text{ for all } x
\end{equation}
On the other hand, we clearly have
\begin{equation}
\phi(\hat{x}(t,y_1,y_2))=0
\end{equation}
it follows that
\begin{equation}
\frac{\partial\hat{\mathcal{H}}_1}{\partial x}(\hat{x}(t,y_1,y_2))=\int_{\mathbb{R}}\frac{\partial \hat{H}_1}{\partial x}(t,\hat{x}(t,y_1,y_2),\hat{u}_1(t,y_1),\hat{u}_2(t,y_2),\hat{p}(t,y_1,y_2),\hat{q}(t,y_1,y_2),\hat{r}(t,y_1,y_2,\zeta)dy_2=a
\end{equation}
Combining this with (\ref{eq2.28}), we get
\begin{align}
& \int_{\mathbb{R}}\int_{\mathbb{R}}[J_1(u_1(.,y_1),\hat{u}_2(.,y_2))-J_1(\widehat{u}_1(.,y_1),\hat{u}_2(.,y_2))]dy_1dy_2\nonumber\\
&\leq\int_0^T\int_{\mathbb{R}}[\hat{\mathcal{H}}_1(x(t,y_1,y_2))-\hat{\mathcal{H}}_1(\hat{x}(t,y_1,y_2))-\frac{\partial\hat{\mathcal{H}}_1}{\partial x}(\hat{x}(t,y_1,y_2))(x(t,y_1,y_2)-\hat{x}(t,y_1,y_2))]dy_1dt\nonumber\\
&\leq 0 \text{ since $\hat{\mathcal{H}}_1$ is concave.}
\end{align}
Hence
\begin{equation}\label{eq2.17}
\sup_{u_1(.,y_1)\in\mathcal{A}_1}\int_{\mathbb{R}}\int_{\mathbb{R}}J_1(u_1(.,y_1),\hat{u}_2(.,y_2))dy_1dy_2\leq \int_{\mathbb{R}}\int_{\mathbb{R}}J_1(\hat{u}_1(.,y_1),\hat{u}_2(.,y_2))dy_1dy_2
\end{equation}
\fproof

\subsection{The case when only one of the players is an insider.}
It is useful also to have a formulation in the partly degenerate case when only one of the players, say player number 1, has inside information. Then the control of player 1 is $\mathbb{H}_1$-adapted as before, while player 2 is $\mathbb{F}$-adapted. In this case we define the \emph{Hamiltonians}
 $ H_i:[0,T]\times\mathbb{R}\times\mathbb{R}\times\mathbb{U}\times\mathbb{R}\times\mathbb{R}\times\mathcal{R} \times \Omega \rightarrow \mathbb{R}$ by
\begin{align}\label{eq2.35}
&H_i(t,x,y_1,u_1,u_2,p,q,r)=H(t,x,y_1,u_1,u_2,p,q,r,\omega)\nonumber\\
&=\mathbb{E}[\delta_{Y_1}(y_1)|\mathcal{F}_t] f_i(t,x,u_1,u_2,y_1)+b(t,x,u_1,u_2,y_1)p \nonumber\\
&+ \sigma(t,x,u_1,u_2,y_1)q+\int_{\mathbb{R}}\gamma(t,x,u_1,u_2,y_1)r(t,\zeta)\nu(d\zeta); i=1,2.
\end{align}
Here, as before, $\mathcal{R}$ denotes the set of all functions $r(.) : \mathbb{R}\rightarrow  \mathbb{R}$
such that the last integral above converges.
For $i=1,2$ we define the \emph{adjoint} processes $p_i(t,y_1),q_i(t,y_1), r_i(t,y_1,\zeta)$ as the solution of the $y_1$-parametrised BSDEs
\begin{equation}\label{eq2.36}
    \left\{
\begin{array}{l}
    dp_i(t,y_1) = -\frac{\partial H_i}{\partial x}(t,y_1)dt + q_i(t,y_1)dB(t)+\int_{\mathbb{R}}r_i(t,y_1,\zeta)\tilde{N}(dt,d\zeta); \quad 0 \leq t\leq T\\
    p_i(T,y_1)  =  g'_i(x(T,y_1),y_1) \mathbb{E}[\delta_{Y_1}(y_1)|\mathcal{F}_T]
\end{array}
    \right.
\end{equation}

Let $J_i(u(.,y_1))$ be defined by
\begin{align}\label{eq2.37}
    J_i(u(.,y_1))&=\mathbb{E}[\int_0^T f_i(t, x(t,y_1),u_1(t,y_1),u_2(t),y_1)\mathbb{E}[\delta_{Y_1}(y_1)|\mathcal{F}_t]dt \nonumber\\
    &+g_i(x(T,y_1),y_1)\mathbb{E}[\delta_{Y_1}(y_1)|\mathcal{F}_T]]
\end{align}
Then we see that
\begin{align}\label{eq2.38}
J_i(u_1,u_2)= \int_{\mathbb{R}} \int_{\mathbb{R}} J_i(u(.,y_1) )dy_1.
\end{align}

\begin{theorem}{[Sufficient maximum principle with only one insider]}\label{sufficient theorem}\\
Suppose $Y_2=0$, i.e. player number 2 has no inside information. Let $(\hat{u_1},\hat{u_2}) \in \mathcal{A}_1\times\mathcal{A}_2$ with associated solution $\hat{x}(t,y_1),\hat{p}_i(t,y_1),\hat{q}_i(t,y_1),\hat{r}_i(t,y_1,\zeta)$ of \eqref{eq8} and \eqref{eq12}; i=1,2. Assume that the following hold:
\begin{enumerate}
 \item $x \rightarrow g_i(x)$ is concave; $i=1,2$
 \item
 The functions
 \begin{equation}
 \hat{\mathcal{H}}_1(x)= \sup_{u_1\in\mathcal{A}_1}H_1(t,x,y_1,u_1,\hat{u_2}(t),\widehat{p}_1(t,y_1),\widehat{q}_1(t,y_1),\hat{r}_1(t,y_1,\cdot))
 \end{equation}
 and
 \begin{equation}
 \hat{\mathcal{H}}_2(x)= \sup_{u_2\in\mathcal{A}_2}\int_{\mathbb{R}}H_2(t,x,y_1,\hat{u}_1(t,y_1),u_2,\widehat{p}_2(t,y_1),\widehat{q}_2(t,y_1),\hat{r}_2(t,y_1,\cdot))dy_1
 \end{equation}
 are concave for all $t.$

\item
 \begin{align}
 \sup_{u_1\in\mathbf{A}_1}&H_1\big(t,\widehat{x}(t,y_1),u_1,\hat{u_2}(t),\widehat{p}_1(t,y_1),\widehat{q}_1(t,y_1),\hat{r}_1(t,y_1,\cdot)\big)\nonumber\\
      &=H_1\big(t,\widehat{x}(t,y_1),\widehat{u}_1(t,y_1),\hat{u_2}(t),\widehat{p}_1(t,y_1),\widehat{q}_1(t,y_1),\hat{r}_1(t,y_1,\cdot)\big)\nonumber\\
      & \text{ for all } t.
\end{align}
\item
 \begin{align}
 \sup_{u_2\in\mathbf{A}_2}&\int_{\mathbb{R}}H_2\big(t,\widehat{x}(t,y_1),\hat{u}_1(t,y_1),u_2,\widehat{p}_2(t,y_1),\widehat{q}_2(t,y_1),\hat{r}_2(t,y_1,\cdot)\big)dy_1\nonumber\\
  &=\int_{\mathbb{R}}H_2\big(t,\widehat{x}(t,y_1),\widehat{u}_1(t,y_1),\hat{u_2}(t),\widehat{p}_2(t,y_1),\widehat{q}_2(t,y_1),\hat{r}_2(t,y_1,\cdot)\big) dy_1\nonumber\\
  &\text{ for all } t.
\end{align}

 \end{enumerate}
Then $(u^*_1(.,y_1),u^*_2(.)) := (\widehat{u}_1(.,y_1),\widehat{u}_2(.))$ is a Nash equilibrium for the problem \eqref{eq2.15}-\eqref{eq2.16}.
\end{theorem}
\dproof
The proof is similar to the proof of Theorem 2.2 and is omitted.
\fproof
\section{A necessary maximum principle}
We proceed to establish a corresponding necessary maximum principle. For this, we do not need concavity conditions, but in stead we need the following assumptions about the set of admissible control values:\\
\begin{itemize}
\item
$A_1$. For all $t_0\in [0,T],y_i \in \mathbb{R}$ and all bounded $\mathcal{F}_{t_0}$-measurable random variables $\alpha_i(y_i,\omega)$, the control
$\theta_i(t,y_i, \omega) := \mathbf{1}_{[t_0,T ]}(t)\alpha_i(y_i,\omega)$ belongs to $\mathcal{A}_i$ for $i=1,2.$\\
\item
$A_2$. For all $u_i; \beta_0^i \in\mathcal{A}_i$ with $\beta_0^i(t,y_i) \leq K < \infty$ for all $t,y_i$  define
\begin{equation}\label{delta}
    \delta_i(t,y_i)=\frac{1}{2K}dist((u_i(t,y_i),\partial\mathbb{A}_i)\wedge1 > 0
\end{equation}
and put
\begin{equation}\label{beta(t,y)}
    \beta_i(t,y_i)=\delta_i(t,y_i)\beta_0^i(t,y_i).
\end{equation}
Then there exists $\delta > 0$ such that the control
\begin{equation*}
    \widetilde{u}_i(t,y_i)=u_i(t,y_i) + a\beta_i(t,y_i) ; \quad t \in [0,T]
\end{equation*}
belongs to $\mathcal{A}_i$ for all $a \in (-\delta, \delta)$ for $i=1,2$.\\
\item
$A3$. For all $\beta_i$ as in (\ref{beta(t,y)}) the derivative processes
\begin{equation*}
    \chi_1(t,y_1,y_2):=\frac{d}{da}x^{(u_1+a\beta_1,u_2)}(t,y_1,y_2)|_{a=0}
\end{equation*}
and
\begin{equation*}
    \chi_2(t,y_1,y_2):=\frac{d}{da}x^{(u_1,u_2+a\beta_2)}(t,y_1,y_2)|_{a=0}
\end{equation*}

exists, and belong to $\mathbf{L}^{2}(\lambda\times \mathbf{P})$ and
\begin{equation}\label{d chi}
    \left\{
\begin{array}{l}
    d\chi_1(t,y_1,y_2) = [\frac{\partial b}{\partial x}(t,y_1,y_2)\chi_1(t,y_1,y_2)+\frac{\partial b}{\partial u_1}(t,y)\beta_1(t,y_1)]dt\\
    +[\frac{\partial \sigma}{\partial x}(t,y_1,y_2)\chi_1(t,y)+\frac{\partial\sigma}{\partial u_1}(t,y_1,y_2)\beta_1(t,y_1)]d B(t) \\
    +\int_{\mathbb{R}}[\frac{\partial \gamma}{\partial x}(t,y_1,y_2,\zeta)\chi_1(t,y_1,y_2)+\frac{\partial \gamma}{\partial u_1}(t,y_1,y_2,\zeta)\beta_1(t,y_1)]\tilde{N}(dt,d\zeta)\\
    \chi_1(0,y_1,y_2)  = \frac{d}{da}x^{(u_1+a\beta_1,u_2)}(0,y_1,y_2)|_{a=0} = 0.
   \end{array}
    \right.
\end{equation}
and
\begin{equation}\label{d chi2}
    \left\{
\begin{array}{l}
    d\chi_2(t,y) = [\frac{\partial b}{\partial x}(t,y_1,y_2)\chi_2(t,y_1,y_2)+\frac{\partial b}{\partial u_2}(t,y_1,y_2)\beta_2(t,y_2)]dt\\
    +[\frac{\partial \sigma}{\partial x}(t,y_1,y_2)\chi_2(t,y_1,y_2)+\frac{\partial\sigma}{\partial u_2}(t,y_1,y_2)\beta_2(t,y_2)]d B(t) \\
    +\int_{\mathbb{R}}[\frac{\partial \gamma}{\partial x}(t,y_1,y_2,\zeta)\chi_2(t,y_1,y_2)+\frac{\partial \gamma}{\partial u_2}(t,y_1,y_2,\zeta)\beta_2(t,y_2)]\tilde{N}(dt,d\zeta)\\
    \chi(0,y_1,y_2)  = \frac{d}{da}x^{(u_1,u_2+a\beta_2)}(0,y_1,y_2)|_{a=0} = 0.
   \end{array}
    \right.
\end{equation}
\end{itemize}

\begin{theorem}{[Necessary maximum principle]} \label{necessary}\\
Let $(u_1,u_2) \in \mathcal{A}_1\times\mathcal{A}_2$. Then the following are equivalent:
\begin{enumerate}
\item $\frac{d}{da}\int_{\mathbb{R}}\int_{\mathbb{R}}J_1(u_1+a\beta_1,u_2)|_{a=0}dy_1dy_2=\frac{d}{da}\int_{\mathbb{R}}\int_{\mathbb{R}}J_2(u_1,u_2+a\beta_2)|_{a=0}dy_1dy_2 =0$ for all bounded $\beta_i \in \mathcal{A}_i$ of the form (\ref{beta(t,y)}).

\item
\begin{align}
&[\int_{\mathbb{R}}\frac{\partial H_1}{\partial v_1}(t,x(t,y_1,y_2),v_1,u_2(t,y_2),p_1(t,y_1,y_2),q_1(t,y_1,y_2),r_1(t,y_1,y_2,.))dy_2]_{v_1=u_1(t,y_1)}\nonumber\\
&=[\int_{\mathbb{R}}\frac{\partial H_2}{\partial v_2}(t,x(t,y_1,y_2),u_1(t,y_1),v_2,,p_2(t,y_1,y_2),q_2(t,y_1,y_2),r_2(t,y_1,y_2,.))dy_1]_{v_2=u_2(t,y_2)}\nonumber\\
& =0\quad \forall t\in[0,T]
 \end{align}
 \end{enumerate}
\end{theorem}

\dproof
By considering an increasing sequence of stopping times $\tau_n$ converging to $T$, we may assume that all local integrals appearing in the computations below are martingales and have expectation $0$. See \cite{OS2}.\\
We can write
$$\frac{d}{da}\int_{\mathbb{R}}\int_{\mathbb{R}}J_1(u_1+a\beta_1,u_2)|_{a=0}dy_1dy_2=I_1+I_2$$\\
where
\begin{align}
I_1&=\frac{d}{da}\int_{\mathbb{R}}\int_{\mathbb{R}}\mathbb{E}[\int_0^Tf_1(t,x^{u_1+a\beta_1}(t,y_1,y_2),u_1(t,y_1)+a\beta_1(t,y_1),u_2(t,y_2),y_1,y_2)\nonumber\\
&\mathbb{E}[\delta_{Y_1,Y_2}(y_1,y_2)|\mathcal{F}_t]dt]|_{a=0}dy_1dy_2
\end{align}
and
$$I_2=\frac{d}{da}\int_{\mathbb{R}}\int_{\mathbb{R}}\mathbb{E}[g_1(x^{(u_1+a\beta_1,u_2)}(T,y_1,y_2),y_1,y_2)\mathbb{E}[\delta_{Y_1,Y_2}(y_1,y_2)|\mathcal{F}_T]]|_{a=0}dy_1dy_2.$$
By our assumptions on $f_1$ and $g_1$ and by \eqref{eq12} we have
\begin{equation}\label{iii1}
    I_1=\int_{\mathbb{R}}\int_{\mathbb{R}}\mathbb{E}[\int_0^T\{\frac{\partial f_1}{\partial x}(t,y_1,y_2)\chi_1(t,y_1,y_2)+\frac{\partial f_1}{\partial u_1}(t,y_1,y_2)\beta_1(t,y_1)\}\mathbb{E}[\delta_{Y_1,Y_2}(y_1,y_2)|\mathcal{F}_t]dt]dy_1dy_2
\end{equation}
\begin{align}\label{iii2}
    I_2&=\int_{\mathbb{R}}\int_{\mathbb{R}}\mathbb{E}[g'_1(x(T,y_1,y_2),y_1,y_2)\chi_1(T,y_1,y_2)\mathbb{E}[\delta_{Y_1,Y_2}(y_1,y_2)|\mathcal{F}_T]]dy_1dy_2\nonumber\\
    &=\int_{\mathbb{R}}\int_{\mathbb{R}}\mathbb{E}[p_1(T,y_1,y_2)\chi_1(T,y_1,y_2)]dy_1dy_2
\end{align}
By the It\^{o} formula
\begin{eqnarray}\label{iii22}
   I_2&=& \int_{\mathbb{R}}\int_{\mathbb{R}}\mathbb{E}[p_1(T,y_1,y_2)\chi_1(T,y_1,y_2)]dy_1dy_2\nonumber\\
   &=&\int_{\mathbb{R}}\int_{\mathbb{R}}\mathbb{E}[\int_0^Tp_1(t,y_1,y_2)d\chi_1(t,y_1,y_2)+\int_0^T\chi_1(t,y_1,y_2)dp_1(t,y_1,y_2)+\int_0^Td[\chi_1,p_1](t,y_1,y_2)]dy_1dy_2 \nonumber\\
   &=&\int_{\mathbb{R}}\int_{\mathbb{R}}\mathbb{E}[\int_0^Tp_1(t,y_1,y_2)\{\frac{\partial b}{\partial x}(t,y_1,y_2)\chi_1(t,y_1,y_2)+\frac{\partial b}{\partial u_1}(t,y_1,y_2)\beta_1(t,y_1)\}dt\\ \nonumber
   &+&\int_0^Tp_1(t,y_1,y_2)\{\frac{\partial \sigma}{\partial x}(t,y_1,y_2)\chi_1(t,y_1,y_2)+\frac{\partial\sigma}{\partial u_1}(t,y_1,y_2)\beta_1(t,y_1)\}dB(t)\nonumber \\
   &+&\int_0^T\int_{\mathbb{R}}p_1(t,y_1,y_2)\{\frac{\partial \gamma}{\partial x}(t,y_1,y_2,\zeta)\chi_1(t,y_1,y_2)+\frac{\partial \gamma}{\partial u_1}(t,y_1,y_2,\zeta)\beta_1(t,y_1)\}\tilde{N}(dt,d\zeta)\nonumber\\
   &-&\int_0^T\chi_1(t,y_1,y_2)\frac{\partial H_1}{\partial x}(t,y_1,y_2) dt+\int_0^T\chi_1(t,y_1,y_2)q_1(t,y_1,y_2)dB(t)\nonumber\\
   &+&\int_0^T\int_{\mathbb{R}}\chi_1(t,y_1,y_2)r_1(t,y_1,y_2,\zeta) \tilde{N}(dt,d\zeta)  \\ \nonumber
     &+&\int_0^Tq_1(t,y_1,y_2) \{\frac{\partial \sigma}{\partial x}(t,y)\chi_1(t,y_1,y_2)+\frac{\partial\sigma}{\partial u_1}(t,y_1,y_2)\beta_1(t,y_1)\}dt \nonumber \\
   &+&\int_0^T\int_{\mathbb{R}}\{\frac{\partial \gamma}{\partial x}(t,y_1,y_2,\zeta)\chi_1(t,y_1,y_2)+\frac{\partial \gamma}{\partial u_1}(t,y_1,y_2,\zeta)\beta_1(t,y_1)\}r_1(t,y_1,y_2,\zeta)\nu(\zeta)dt]dy_1dy_2\nonumber\\
  &=& \int_{\mathbb{R}}\int_{\mathbb{R}}\mathbb{E}[\int_0^T\chi_1(t,y_1,y_2)\{p_1(t,y_1,y_2)\frac{\partial b}{\partial x}(t,y_1,y_2)+q_1(t,y_1,y_2)\frac{\partial \sigma}{\partial x}(t,y_1,y_2)-\frac{\partial H_1}{\partial x}(t,y_1,y_2)\nonumber\\
  &+&\int_{\mathbb{R}}\frac{\partial \gamma}{\partial x}(t,y_1,y_2,\zeta)r_1(t,y_1,y_2,\zeta)\nu(d\zeta)\}dt \nonumber \\
   &+& \int_0^T\beta_1(t,y_1)\{p_1(t,y_1,y_2)\frac{\partial b}{\partial u_1}(t,y_1,y_2)+q_1(t,y_1,y_2)\frac{\partial\sigma}{\partial u_1}(t,y_1,y_2)\nonumber\\
   &+&\int_{\mathbb{R}}\frac{\partial \gamma}{\partial u_1}(t,y_1,y_2,\zeta)r_1(t,y_1,y_2,\zeta)\nu(d\zeta)\}dt]dy_1dy_2\nonumber\\
&=&\int_{\mathbb{R}}\int_{\mathbb{R}}\mathbb{E}[-\int_0^T\chi_1(t,y_1,y_2)\frac{\partial f_1}{\partial x}\mathbb{E}[\delta_{Y_1,Y_2}(y_1,y_2)|\mathcal{F}_t]dt\nonumber\\
&+&\int_0^T\{\frac{\partial H_1}{\partial u_1}(t,y_1,y_2)-\frac{\partial f_1}{\partial u_1}(t,y_1,y_2)\mathbb{E}[\delta_{Y_1,Y_2}(y_1,y_2)|\mathcal{F}_t]\}\beta_1(t,y_1)dt]dy_1dy_2\nonumber\\
&=& -I_1+\int_{\mathbb{R}}\int_{\mathbb{R}}\mathbb{E}[\int_0^T\frac{\partial H_1}{\partial u_1}(t,y_1,y_2)\beta_1(t,y_1)dt]dy_1dy_2.
\end{eqnarray}

Summing (\ref{iii1}) and (\ref{iii22}) we get
\begin{equation*}
    \frac{d}{da}\int_{\mathbb{R}}\int_{\mathbb{R}}J_1(u_1+a\beta_1,u_2)|_{a=0}dy_1dy_2=I_1+I_2=\int_{\mathbb{R}}\int_{\mathbb{R}}\mathbb{E}[\int_0^T\frac{\partial H_1}{\partial u_1}(t,y_1,y_2)\beta_1(t,y_1)dt]dy_1dy_2.
\end{equation*}
we conclude that
\begin{equation*}
    \frac{d}{da}\int_{\mathbb{R}}\int_{\mathbb{R}}J_1(u_1+a\beta_1,u_2)|_{a=0}dy_1dy_2=0
\end{equation*}
if and only if
$$\int_{\mathbb{R}}\int_{\mathbb{R}}\mathbb{E}[\int_0^T\frac{\partial H_1}{\partial u_1}(t,y_1,y_2)\beta_1(t,y_1)dt]dy_1dy_2=0$$
 for all bounded $\beta_1\in\mathcal{A}_1$ of the form (\ref{beta(t,y)}).\\
Changing the order of integration we can write this as follows:
\begin{equation}
\mathbb{E}[\int_0^T\int_{\mathbb{R}}F_1(t,y_1)\beta_1(t,y_1)dtdy_1]=0,\quad \forall \beta_1\in\mathcal{A}_1
\end{equation}
where
\begin{equation}
F_1(t,y_1):=\int_{\mathbb{R}}\frac{\partial H_1}{\partial u_1}(t,y_1,y_2)dy_2.
\end{equation}
\noindent In particular, applying this to $\beta_1(t,y_1) = \theta_1(t,y_1)$ as in $A1$, we get that this is again equivalent to
\begin{equation*}
\mathbb{E}[F_1(t,y_1)|\mathcal{F}_t]=0, \quad \forall t,y_1.
\end{equation*}
Since $F_1(t,y_1)$ is already $\mathcal{F}_t$-adapted, we have
\begin{equation*}
\mathbb{E}[F_1(t,y_1)|\mathcal{F}_t]=F_1(t,y_1), \quad\forall t,y_1.
\end{equation*}
So we deduce that
\begin{equation*}
F_1(t,y_1) = \int_{\mathbb{R}}\frac{\partial H_1}{\partial u_1}(t,y_1,y_2)dy_2=0, \quad\forall t,y_1.
\end{equation*}
A similar argument gives that

\begin{equation}
\frac{d}{da}\int_{\mathbb{R}}\int_{\mathbb{R}}J_2(u_1,u_2+a\beta_2)|_{a=0}dy_1dy_2 =0 \text{ for all bounded } \beta_2 \in \mathcal{A}_2
\end{equation}
is equivalent to
\begin{equation}
\int_{\mathbb{R}}\frac{\partial H_2}{\partial u_2}(t,y_1,y_2)dy_1=0, \quad\forall t,y_2,
\end{equation}
where
\begin{align}
&\frac{\partial H_2}{\partial u_2}(t,y_1,y_2)\nonumber\\
 &=\frac{\partial H_2}{\partial v_2}(t,x(t,y_1,y_2),u_1(t,y_1),v_2,p_2(t,y_1,y_2),q_2(t,y_1,y_2),r_2(t,y_1,y_2,.))_{v_2=u_2(t,y_2)}.
 \end{align}
\fproof
\section{The zero-sum game case}
In the zero-sum case we have
\begin{equation}
\int_{\mathbb{R}}\int_{\mathbb{R}}J_1(u_1(.,y_1),u_2(.,y_2))+J_2(u_1(.,y_1),u_2(.,y_2))dy_1dy_2=0.
\end{equation}
Then the Nash equilibrium $(\hat{u}_1(.,y_1),\hat{u}_2(.,y_2))\in\mathcal{A}_1\times\mathcal{A}_2$ satisfying (\ref{eq2.15})-(\ref{eq2.16}) becomes a saddle point for
\begin{equation}
\int_{\mathbb{R}}\int_{\mathbb{R}}J(u_1(.,y_1),u_2(.,y_2))dy_1dy_2:=\int_{\mathbb{R}}\int_{\mathbb{R}}J_1(u_1(.,y_1),u_2(.,y_2))dy_1dy_2.
\end{equation}
To see this, note that (\ref{eq2.15})-(\ref{eq2.16}) imply that
\begin{align}
\int_{\mathbb{R}^2}J_1(u_1(.,y_1),\hat{u}_2(.,y_2))dy_1dy_2\leq \int_{\mathbb{R}^2}J_1(\hat{u}_1(.,y_1),\hat{u}_2(.,y_2))dy_1dy_2&=-\int_{\mathbb{R}^2}J_2(\hat{u}_1(.,y_1),\hat{u}_2(.,y_2))dy_1dy_2\nonumber\\
&\leq -\int_{\mathbb{R}^2}J_2(\hat{u}_1(.,y_1),u_2(.,y_2))dy_1dy_2
\end{align}
and hence
\begin{equation}
\int_{\mathbb{R}^2}J(u_1(.,y_1),\hat{u}_2(.,y_2))dy_1dy_2\leq \int_{\mathbb{R}^2}J(\hat{u}_1(.,y_1),\hat{u}_2(.,y_2))dy_1dy_2\leq \int_{\mathbb{R}^2}J(\hat{u}_1(.,y_1),u_2(.,y_2))dy_1dy_2
\end{equation}
for all $u_1, u_2$.
From this we deduce that
\begin{align}
&\inf_{u_2\in\mathcal{A}_2}\sup_{u_1\in\mathcal{A}_1}\int_{\mathbb{R}^2}J(u_1(.,y_1),u_2(.,y_2))dy_1dy_2\leq
\sup_{u_1\in\mathcal{A}_1}\int_{\mathbb{R}^2}J(u_1(.,y_1),\hat{u}_2(.,y_2))dy_1dy_2\nonumber\\
&\leq \int_{\mathbb{R}^2}J(\hat{u}_1(.,y_1),\hat{u}_2(.,y_2))dy_1dy_2
\leq \inf_{u_2\in\mathcal{A}_2}\int_{\mathbb{R}^2}J(\hat{u}_1(.,y_1),u_2(.,y_2))dy_1dy_2\nonumber\\
&\leq \sup_{u_1\in\mathcal{A}_1}\inf_{u_2\in\mathcal{A}_2}\int_{\mathbb{R}^2}J(u_1(.,y_1),u_2(.,y_2))dy_1dy_2.
\end{align}
Since we always have $\inf\sup\geq\sup\inf$, we conclude that
\begin{align}
&\inf_{u_2\in\mathcal{A}_2}\sup_{u_1\in\mathcal{A}_1}\int_{\mathbb{R}^2}J(u_1(.,y_1),u_2(.,y_2))dy_1dy_2=\sup_{u_1\in\mathcal{A}_1}\int_{\mathbb{R}^2}J(u_1(.,y_1),\hat{u}_2(.,y_2))dy_1dy_2\nonumber\\
&= \int_{\mathbb{R}^2}J(\hat{u}_1(.,y_1),\hat{u}_2(.,y_2))dy_1dy_2
=\inf_{u_2\in\mathcal{A}_2}\int_{\mathbb{R}^2}J(\hat{u}_1(.,y_1),u_2(.,y_2))dy_1dy_2\nonumber\\
&=\sup_{u_1\in\mathcal{A}_1}\inf_{u_2\in\mathcal{A}_2}\int_{\mathbb{R}^2}J(u_1(.,y_1),u_2(.,y_2))dy_1dy_2
\end{align}
i.e $(\hat{u}_1(.,y_1),\hat{u}_2(.,y_2))\in\mathcal{A}_1\times\mathcal{A}_2$ is a saddle point for $\int_{\mathbb{R}}\int_{\mathbb{R}}J(u_1(.,y_1),u_2(.,y_2))dy_1dy_2$.
Hence we want to find $(\hat{u}_1(.,y_1),\hat{u}_2(.,y_2))\in\mathcal{A}_1\times\mathcal{A}_2$ such that
\begin{align}
\sup_{u_1\in\mathcal{A}_1}\inf_{u_2\in\mathcal{A}_2}\int_{\mathbb{R}^2}J(u_1(.,y_1),u_2(.,y_2))dy_1dy_2&=\inf_{u_2\in\mathcal{A}_2}\sup_{u_1\in\mathcal{A}_1}\int_{\mathbb{R}^2}J(u_1(.,y_1),u_2(.,y_2))dy_1dy_2\nonumber\\
&=\int_{\mathbb{R}}\int_{\mathbb{R}}J(\hat{u}_1(.,y_1),\hat{u}_2(.,y_2))dy_1dy_2
\end{align}
where
\begin{align}\label{J(u)}
    \int_{\mathbb{R}}\int_{\mathbb{R}}J(u(.,y_1,y_2))dy_1dy_2&=\int_{\mathbb{R}}\int_{\mathbb{R}}\mathbb{E}\big[\int_0^T f(t, x(t,y_1,y_2),u_1(t,y_1),u_2(t,y_2),y_1,y_2)\mathbb{E}[\delta_{Y_1,Y_2}(y_1,y_2)|\mathcal{F}_t]dt \nonumber\\
    &+g(x(T,y_1,y_2),y_1,y_2)\mathbb{E}[\delta_{Y_1,Y_2}(y_1,y_2)|\mathcal{F}_T]\big]dy_1dy_2
\end{align}
\subsection{Situation 1: Both players are still maximising their own performance functional }

Choose $g_1=g=-g_2$ and $f_1=f=-f_2$. Then by (\ref{eq11}) the Hamiltonians are:

\begin{align}\label{eq4.11}
&H_1(t,x,y_1,y_2,u_1,u_2,p,q,r)=H_1(t,x,y_1,y_2,u_1,u_2,p,q,r,\omega)\nonumber\\
&=\mathbb{E}[\delta_{Y_1,Y_2}(y_1,y_2)|\mathcal{F}_t] f(t,x,u_1,u_2,y_1,y_2)+b(t,x,u_1,u_2,y_1,y_2)p \nonumber\\
&+ \sigma(t,x,u_1,u_2,y_1,y_2)q+\int_{\mathbb{R}}\gamma(t,x,u_1,u_2,y_1,y_2)r(t,\zeta)\nu(d\zeta)
\end{align}
and
\begin{align}\label{eq4.12}
&H_2(t,x,y_1,y_2,u_1,u_2,p,q,r)=H_2(t,x,y_1,y_2,u_1,u_2,p,q,r,\omega)\nonumber\\
&=-\mathbb{E}[\delta_{Y_1,Y_2}(y_1,y_2)|\mathcal{F}_t] f(t,x,u_1,u_2,y_1,y_2)+b(t,x,u_1,u_2,y_1,y_2)p \nonumber\\
&+ \sigma(t,x,u_1,u_2,y_1,y_2)q+\int_{\mathbb{R}}\gamma(t,x,u_1,u_2,y_1,y_2)r(t,\zeta)\nu(d\zeta)
\end{align}
Let $p_i, q_i,r_i$,  $i = 1, 2$ be as in (\ref{eq12}).

Let us now state the necessary maximum principle for the zero sum game problem:
\begin{theorem}{[Necessary maximum principle for zero-sum games]}\\
Assume the conditions of Theorem \ref{necessary} hold.
Then the following are equivalent:
\begin{enumerate}
\item $\frac{d}{da}\int_{\mathbb{R}}\int_{\mathbb{R}}J(u_1+a\beta_1,u_2)|_{a=0}dy_1dy_2=\frac{d}{da}\int_{\mathbb{R}}\int_{\mathbb{R}}J(u_1,u_2+a\beta_2)|_{a=0}dy_1dy_2 =0$ for all bounded $\beta_i \in \mathcal{A}_i$ of the form (\ref{beta(t,y)}).
\item
\begin{align}
&[\int_{\mathbb{R}}\frac{\partial H_1}{\partial v_1}(t,x(t,y_1,y_2),v_1,u_2(t,y_2),p_1(t,y_1,y_2),q_1(t,y_1,y_2),r_1(t,y_1,y_2,.))dy_2]_{v_1=u_1(t,y_1)}\nonumber\\
&=[\int_{\mathbb{R}}\frac{\partial H_2}{\partial v_2}(t,x(t,y_1,y_2),u_1(t,y_1),v_2,,p_2(t,y_1,y_2),q_2(t,y_1,y_2),r_2(t,y_1,y_2,.))dy_1]_{v_2=u_2(t,y_2)}\nonumber\\
& =0\quad \forall t\in[0,T],y_1,y_2.
 \end{align}
 \end{enumerate}
\end{theorem}

\dproof
This is a direct consequence of Theorem \ref{necessary}
\fproof

\subsection{Situation 2: One of the players is maximising and the other is minimising the performance functional}
 Let us now look at the problem with one performance functional common to both players, but where one of the players is maximising and the other is minimising it. Then we get just one Hamiltonian and just one BSDE, which is simpler to deal with.

In this case the  Hamiltonian $H$, is given by:
\begin{align}\label{eq11}
&H(t,x,y_1,y_2,u_1,u_2,p,q,r)=H(t,x,y_1,y_2,u_1,u_2,p,q,r,\omega)\nonumber\\
&=\mathbb{E}[\delta_{Y_1,Y_2}(y_1,y_2)|\mathcal{F}_t] f(t,x,u_1,u_2,y_1,y_2)+b(t,x,u_1,u_2,y_1,y_2)p \nonumber\\
&+ \sigma(t,x,u_1,u_2,y_1,y_2)q+\int_{\mathbb{R}}\gamma(t,x,u_1,u_2,y_1,y_2)r(t,\zeta)\nu(d\zeta)
\end{align}
Moreover, there is only one triple $(p,q,r)$ of adjoint processes, given by the BSDE
\begin{equation}\label{eq4.12}
    \left\{
\begin{array}{l}
    dp(t,y_1,y_2) = -\frac{\partial H}{\partial x}(t,y_1,y_2)dt +q(t,y_1,y_2)dB(t)+\int_{\mathbb{R}}r(t,y_1,y_2,\zeta)\tilde{N}(dt,d\zeta); \quad 0 \leq t\leq T\\
    p(T,y)  = g'(x(T,y_1,y_2),y_1,y_2) \mathbb{E}[\delta_{Y_1,Y_2}(y_1,y_2)|\mathcal{F}_T]
\end{array}
    \right.
\end{equation}
We can now state the corresponding sufficient maximum principle for the zero-sum game:
\begin{theorem}(Sufficient maximum principle for the zero-sum game)

Let $(\hat{u_1},\hat{u_2}) \in \mathcal{A}_1\times\mathcal{A}_2$ with associated solution $\hat{x}(t,y_1,y_2),\hat{p}(t,y_1,y_2),\hat{q}(t,y_1,y_2),\hat{r}(t,y_1,y_2,\zeta)$ of \eqref{eq8} and \eqref{eq4.12}. Assume that the following holds:

\begin{enumerate}

\item the function $x \rightarrow g(x)$ is affine
\item
\begin{align}
&\sup_{u_1\in\mathbf{A}_1}\int_{\mathbb{R}}H\big(t,\widehat{x}(t,y_1,y_2),u_1,\hat{u_2}(t,y_2),\widehat{p}(t,y_1,y_2),\widehat{q}(t,y_1,y_2),\hat{r}(t,y_1,y_2,\cdot)\big)dy_2\nonumber\\
&=\int_{\mathbb{R}}H\big(t,\widehat{x}(t,y_1,y_2),\widehat{u}_1(t,y_1),\hat{u_2}(t,y_2),\widehat{p}(t,y_1,y_2),\widehat{q}(t,y_1,y_2),\hat{r}(t,y_1,y_2,\cdot)\big)dy_2\nonumber\\
&\text{ for all } t,y_1.
\end{align}
\begin{align}
&\inf_{u_2\in\mathbf{A}_2}\int_{\mathbb{R}}H\big(t,\widehat{x}(t,y_1,y_2),\hat{u}_1(t,y_1),u_2,\widehat{p}(t,y_1,y_2),\widehat{q}(t,y_1,y_2),\hat{r}(t,y_1,y_2,\cdot)\big)dy_1\nonumber\\
&=\int_{\mathbb{R}}H\big(t,\widehat{x}(t,y_1,y_2),\hat{u}_1(t,y_1),\hat{u}_2(t,y_2),\widehat{p}(t,y_1,y_2),\widehat{q}(t,y_1,y_2),\hat{r}(t,y_1,y_2,\cdot)\big)dy_1\nonumber\\
&\text{ for all } t,y_2.
\end{align}

 \item The function

\begin{equation}
\hat{\mathcal{H}}(x)=\sup_{u_1\in\mathcal{A}_1}\int_{\mathbb{R}} H(t,x,y_1,y_2,u_1,\hat{u_2}(t,y_2),\widehat{p}(t,y_1,y_2),\widehat{q}(t,y_1,y_2),\hat{r}(t,y_1,y_2,\cdot))dy_2
\end{equation}
is concave for all $t,y_1,$\\
and the function
\begin{equation}
 \underline{\mathcal{H}}(x)= \inf_{u_2\in\mathcal{A}_2}\int_{\mathbb{R}}H(t,x,y_1,y_2,\hat{u}_1(t,y_1),u_2,\widehat{p}(t,y_1,y_2),\widehat{q}(t,y_1,y_2),\hat{r}(t,y_1,y_2,\cdot))dy_1
\end{equation}
is convex for all $t,y_2.$

 \end{enumerate}
 Then
 $\hat{u}(t,y_1,y_2)=(\hat{u}_1(t,y_1),\hat{u}_2(t,y_2))$ is a saddle point for $J(u_1,u_2)$.
\end{theorem}

Let us now state the necessary maximum principle for the zero sum game problem:
\begin{theorem}{[Necessary maximum principle for zero-sum games]}\\
Assume the conditions of Theorem \ref{necessary} hold.
Then the following are equivalent:
\begin{enumerate}
\item $\frac{d}{da}\int_{\mathbb{R}}\int_{\mathbb{R}}J(u_1+a\beta_1,u_2)|_{a=0}dy_1dy_2=\frac{d}{da}\int_{\mathbb{R}}\int_{\mathbb{R}}J(u_1,u_2+a\beta_2)|_{a=0}dy_1dy_2 =0$ for all bounded $\beta_i \in \mathcal{A}_i$ of the form (\ref{beta(t,y)}).

\item
\begin{align}
&[\int_{\mathbb{R}}\frac{\partial H}{\partial v_1}(t,x(t,y_1,y_2),v_1,u_2(t,y_2),p_1(t,y_1,y_2),q_1(t,y_1,y_2),r_1(t,y_1,y_2,.))dy_2]_{v_1=u_1(t,y_1)}\nonumber\\
&=[\int_{\mathbb{R}}\frac{\partial H}{\partial v_2}(t,x(t,y_1,y_2),u_1(t,y_1),v_2,,p_2(t,y_1,y_2),q_2(t,y_1,y_2),r_2(t,y_1,y_2,.))dy_1]_{v_2=u_2(t,y_2)}\nonumber\\
& =0\quad \forall t\in[0,T],y_1,y_2.
 \end{align}
 \end{enumerate}

\end{theorem}

\section{Applications}
\subsection{Optimal insider consumption under model uncertainty}
Suppose we have a cash flow with consumption, modelled by the process $X(t,Y)=X^{c,\mu}(t,Y)$ defined by:
\begin{equation}
\begin{cases}
dX(t,Y)= (\alpha(t,Y) +\mu(t,Y_2) - c(t,Y_1))X(t,Y)dt + \beta(t,Y)X(t,Y)dB(t) +\int_{\mathbb{R}}\gamma(t,Y,\zeta)X(t,Y)\tilde{N}(dt,d\zeta)\nonumber\\
X(0)=x >0
\end{cases}
\end{equation}
Here $\alpha(t,Y),\beta(t,Y),\gamma(t,Y)$ are given coefficients, while $c(t,Y_1)>0$ is the relative consumption rate chosen by the consumer (player number 1) and $\mu(t,Y_2)$ is a perturbation of the drift term, representing the model uncertainty chosen by the environment (player number 2).
Define the performance functional by
\begin{equation}
J(c,\mu)= \mathbb{E}[\int_0^T\{\log(c(t)X(t))+\frac{1}{2} \mu^2(t)\}dt + \theta \log X(T)]
\end{equation}
where $\theta>0$ is a given constant and $\frac{1}{2} \mu^2(t)$ represents a penalty rate, penalizing $\mu$ for being away from 0.
We assume that $c$ is $\mathbb{H}^1$-adapted, while $\mu$ is $\mathbb{H}^2$-adapted.

We want to find $c^* \in \mathcal{A}_1$ and $\mu^* \in \mathcal{A}_2$ such that

\begin{equation}\label{eq0.2}
\sup_{c \in \mathcal{A}_1} \inf_{\mu \in \mathcal{A}_2} J(c,\mu)=J(c^*,\mu^*).
\end{equation}

As before we rewrite this problem as a classical stochastic differential game with two parameters $y_1,y_2$.
Thus we define, for $y=(y_1,y_2) \in \mathbb{R} \times \mathbb{R}$,
\begin{align}
\begin{cases}
dx(t,y)&= (\alpha(t,y) +\mu(t,y_2) - c(t,y_1))x(t,y)dt + \beta(t,y)x(t,y)dB(t) \\
&+\int_{\mathbb{R}}\gamma(t,y,\zeta)x(t,y)\tilde{N}(dt,d\zeta)\\
x(0,y)&=x >0
\end{cases}
\end{align}
and
\begin{align}
&J(c(.,y_1),\mu(.,y_2))\nonumber\\
&=\mathbb{E}[\int_0^T\{\log(c(t,y_1)x(t,y))+\frac{1}{2} \mu^2(t,y_2)\} \mathbb{E}[\delta_Y(y)|\mathcal{F}_t]dt + \theta \log x(T,y)\mathbb{E}[\delta_Y(y)|\mathcal{F}_T]]
\end{align}

The Hamiltonian for this problem is
\begin{align}
&H(t,x,y,c,\mu,p,q,r)\nonumber\\
&= \{\log(cx)+\frac{1}{2}\mu^2\}\mathbb{E}[\delta_Y(y)|\mathcal{F}_t] +(\alpha(t,y)+\mu-c)xp +\beta(t,y)xq+x\int_{\mathbb{R}}\gamma(t,y,\zeta)r(\zeta) d\nu(\zeta)
\end{align}
and the BSDE for the adjoint processes $p,q,r$ is

\begin{align}\label{eq0.6}
\begin{cases}
dp(t,y)&= -[ \frac{1}{x(t,y)}\mathbb{E}[\delta_Y(y)|\mathcal{F}_t]+(\alpha(t,y)+\mu(t,y_2)-c(t,y_1))p(t,y) \\ &+\beta(t,y)q(t,y)+\int_{\mathbb{R}}\gamma(t,y,\zeta)r(\zeta) d\nu(\zeta)]dt \\
&+ q(t,y)dB(t) + \int_{\mathbb{R}}r(t,y)\tilde{N}(dt,d\zeta); 0 \leq t \leq T\\
p(T,y)&=\frac{\theta}{x(T,y)}\mathbb{E}[\delta_Y(y)|\mathcal{F}_T]
\end{cases}
\end{align}
Define
\begin{equation}\label{eq0.7}
h(t,y)=p(t,y)x(t,y).
\end{equation}
Then by the It\^ o formula we get
\begin{align}\label{eq0.8}
&dh(t,y)=x(t,y)[-\frac{1}{x(t,y)}\mathbb{E}[\delta_Y(y)|\mathcal{F}_t]-(\alpha(t,Y) +\mu(t,Y_2) - c(t,Y_1))p(t,y)-\beta(t,y)q(t,y)\nonumber\\
&-\int_{\mathbb{R}}\gamma(t,y,\zeta)r(t,\zeta)d\nu(\zeta)]dt\nonumber\\
&+p(t,y)(\alpha(t,Y) +\mu(t,Y_2) - c(t,Y_1))x(t,y)dt+p(t,y)\beta(t,y)x(t,y)dB(t)+x(t,y)q(t,y)dB(t)\nonumber\\
&+q(t,y)\beta(t,y)x(t,y)dt\nonumber\\
&+\int_{\mathbb{R}}[(x(t,y)+\gamma(t,y,\zeta)x(t,y))(p(t,y)+r(t,y,\zeta))-p(t,y)x(t,y)-p(t,y)\gamma(t,y,\zeta)x(t,y)-x(t,y)r(t,y,\zeta]d\nu(\zeta)dt\nonumber\\\\
&+\int_{\mathbb{R}}[(x(t,y)+\gamma(t,y,\zeta)x(t,y))(p(t,y)+r(t,y,\zeta))-p(t,y)x(t,y)]\tilde{N}(dt,d\zeta)\nonumber\\
&=dF(t,y) +h(t,y)\beta(t,y)dB(t)+h(t,y)\int_{\mathbb{R}}\gamma(t,y,\zeta)\tilde{N}(dt,d\zeta)),
\end{align}
where
\begin{align}
&dF(t,y)=\nonumber\\
& -\mathbb{E}[\delta_Y(y)|\mathcal{F}_t]dt+ x(t,y)q(t,y)dB(t) + x(t,y)\int_{\mathbb{R}}r(t,y,\zeta)(1+\gamma(t,y,\zeta))\tilde{N}(dt,d\zeta).
\end{align}

To simplify this, we define the process $k(t,y)$ by the equation
\begin{align}\label{eq0.11}
&dk(t,y)=k(t,y)\Big[ b(t,y)dB(t)+\int_{\mathbb{R}}c(t,y,\zeta)\tilde{N}(dt,d\zeta)\Big]
\end{align}
for suitable processes $b, c$ (to be determined).

Then again by the It\^ o formula we get
\begin{align}\label{eq0.12}
d(h(t,y)k(t,y))&=h(t,y)k(t,y)\Big[ b(t,y)dB(t)+\int_{\mathbb{R}}c(t,y,\zeta)\tilde{N}(dt,d\zeta)\Big] \nonumber\\
&+ k(t,y)\Big[ dF(t,y)+h(t,y)\beta(t,y)dB(t) +h(t,y)\int_{\mathbb{R}}\gamma(t,y,\zeta)\tilde{N}(dt,d\zeta)\Big]\nonumber\\
&+(h(t,y)\beta(t,y)+x(t,y)q(t,y))k(t,y)b(t,y)dt\nonumber\\
&+\int_{\mathbb{R}}\Big( h(t,y)\gamma(t,y,\zeta) + x(t,y)r(t,y,\zeta)(1+\gamma(t,y,\zeta)\Big) k(t,y)c(t,y,\zeta)\tilde{N}(dt,d\zeta)\nonumber\\
&+\int_{\mathbb{R}}\Big( h(t,y)\gamma(t,y,\zeta) + x(t,y)r(t,y,\zeta)(1+\gamma(t,y,\zeta)\Big) k(t,y)c(t,y,\zeta)d\nu(\zeta)dt
\end{align}

Define
\begin{align}\label{eq0.13}
&u(t,y):=h(t,y)k(t,y). \nonumber\\
\end{align}
Then the equation above can be written
\begin{align}\label{eq0.14}
du(t,y)&= u(t,y)\Big[ \int_{\mathbb{R}}\gamma(t,y,\zeta)c(t,y,\zeta)d\nu(\zeta) dt\nonumber\\
&+\{ \beta(t,y)+b(t,y)\} dB(t)+\beta(t,y)b(t,y)dt +\int_{\mathbb{R}} \{c(t,y,\zeta) + \gamma(t,y,\zeta)+c(t,y,\zeta)\gamma(t,y,\zeta)\} \tilde{N}(dt,d\zeta)\Big]\nonumber\\
&+k(t,y)\Big[ dF(t,y)+x(t,y)q(t,y)b(t,y)dt+\int_{\mathbb{R}}x(t,y)r(t,y,\zeta)c(t,y,\zeta)(1+\gamma(t,y,\zeta))d\nu(\zeta)dt\nonumber\\
&+\int_{\mathbb{R}}x(t,y)r(t,y,\zeta)c(t,y,\zeta)(1+\gamma(t,y,\zeta))\tilde{N}(dt,d\zeta)\Big]
\end{align}
Choose
\begin{align}\label{eq0.15}
&b(t,y):=-\beta(t,y)\nonumber\\
&c(t,\zeta):= -\frac{\gamma(t,y,\zeta)}{1+\gamma(t,y,\zeta)}
\end{align}

Then \eqref{eq0.14} reduces to
\begin{align}\label{eq0.16}
du(t,y)&=f(t,y)dt+k(t,y)x(t,y)q(t,y)dB(t) \nonumber\\
&+ \int_{\mathbb{R}}\{ x(t,y)r(t,y,\zeta)(1+\gamma(t,y,\zeta))[k(t,y)+k(t,y)c(t,y,\zeta)]\} \tilde{N}(dt,d\zeta),
\end{align}
where
\begin{align}\label{eq0.17}
f(t,y)&=-k(t,y)E[\delta_Y(y)|\mathcal{F}_t]+u(t,y)[\int_{\mathbb{R}}\gamma(t,y,\zeta)c(t,\zeta)d\nu(\zeta)+\beta(t,y)b(t,y)]\nonumber\\
&+k(t,y)x(t,y)q(t,y)b(t,y)+k(t,y)\int_{\mathbb{R}}x(t,y)r(t,y,\zeta)c(t,y,\zeta)(1+\gamma(t,y,\zeta))d\nu(\zeta)
\end{align}
Now define
\begin{align}\label{eq0.18}
& v(t,y):= k(t,y)x(t,y)q(t,y)\nonumber\\
&w(t,y):=k(t,y)x(t,y)r(t,y,\zeta).
\end{align}
Then from \eqref{eq0.12} and \eqref{eq0.15} we get the following BSDE in the unknowns $u, v, w$:
\begin{align}\label{eq0.19}
du(t,y)=&\Big( -k(t,y)\mathbb{E}[\delta_Y(y)|\mathcal{F}_t]-u(t,y)[\int_{\mathbb{R}}\frac{\gamma^2(t,y,\zeta)}{1+\gamma(t,y,\zeta)}d\nu(\zeta)+ \beta^2(t,y)]\\
&-\beta(t,y)v(t,y) - \int_{\mathbb{R}}\gamma(t,y,\zeta) w(t,y,\zeta)d\nu(\zeta)\Big) dt\nonumber\\
&+ v(t,y)dB(t)+ \int_{\mathbb{R}}w(t,y,\zeta)\tilde{N}(dt,d\zeta); \quad 0 \leq t \leq T\nonumber\\
u(T,y)= & \theta k(T,y)\mathbb{E}[\delta_Y(y)|\mathcal{F}_T]\nonumber
\end{align}
This is a linear BSDE which has a unique solution $u(t,y)=p(t,y)x(t,y)k(t,y),v(t,y),w(t,y,\zeta)$. In particular, we may regard
\begin{equation}\label{eq0.20}
p(t,y)x(t,y)= \frac{u(t,y)}{k(t,y)}
\end{equation}
as known.

Maximizing $\int_{\mathbb{R}}Hdy_2$ with respect to c gives the first order equation
\begin{align}
\int_{\mathbb{R}} \{\frac{1}{c(t,y_1)}\mathbb{E}[\delta_Y(y)|\mathcal{F}_t] -x(t,y)p(t,y)\}dy_2=0,
\end{align}
i.e.,
\begin{align}\label{eq0.22}
c(t,y_1)=\hat{c}(t,y_1)= \frac{\int_{\mathbb{R}}\mathbb{E}[\delta_Y(y)|\mathcal{F}_t]dy_2}{\int_{\mathbb{R}}x(t,y)p(t,y)dy_2}.
\end{align}

Minimizing $\int_{\mathbb{R}}Hdy_1$ with respect to $\mu$ gives the first order equation
\begin{align}
\int_{\mathbb{R}}\{\mu(t,y_2)\mathbb{E}[\delta_Y(y)|\mathcal{F}_t]+x(t,y)p(t,y)\}dy_1=0,
\end{align}
i.e.
\begin{align}\label{eq0.24}
\mu(t,y_2)=\hat{\mu}(t,y_2)= \frac{\int_{\mathbb{R}}x(t,y)p(t,y)dy_1}{\int_{\mathbb{R}}\mathbb{E}[\delta_Y(y)|\mathcal{F}_t]dy_1}.
\end{align}

We can now verify that $\hat{c},\hat{\mu}$ satisfies all the conditions of the sufficient maximum principle, and hence we conclude the following:

\begin{theorem}[Optimal consumption for an insider under model uncertainty]
The solution $(c^*,\mu^*)$ of the stochastic differential game \eqref{eq0.2} is given by
\begin{align}\label{eq0.25}
c^*(t,Y_1)= \frac{\int_{\mathbb{R}}\mathbb{E}[\delta_Y(y)|\mathcal{F}_t]dy_2|_{y_1=Y_1}}{\int_{\mathbb{R}}x(t,y)p(t,y)dy_2 |_{y_1=Y_1}}.
\end{align}
and
\begin{align}\label{eq0.26}
\mu^*(t,Y_2)= \frac{\int_{\mathbb{R}}x(t,y)p(t,y)dy_1|_{y_2=Y_2}}{\int_{\mathbb{R}}\mathbb{E}[\delta_Y(y)|\mathcal{F}_t]dy_1|_{y_2=Y_2}},
\end{align}

where $h(t,y)=x(t,y)p(t,y)$ is given by \eqref{eq0.19}-\eqref{eq0.20}.

\end{theorem}

\subsection{Optimal insider portfolio under model uncertainty}
Consider a financial market with two investment possibilities:
\begin{itemize}
\item
(i) A risk free investment possibility with unit price $S_0(t)=1$ for all $t \in [0,T]$
\item
(ii) A risky investment, where the unit price $S(t)=S(t,Y)$ is modelled by the (forward) SDE
\begin{equation}\label{eq00.1}
dS(t,Y) = S(t,Y) [(\alpha(t,Y)+\mu(t))dt + \beta(t,Y) dB(t)]; S(0) >0.
\end{equation}
\end{itemize}
Here $\alpha(t,Y),\beta(t,Y)$ are given $\mathbb{H}$-adapted coefficients, while $\mu(t)$ is a perturbation of the drift term, representing the model uncertainty chosen by the environment (player number 2).

Suppose the wealth process $X(t,Y)=X^{\pi,\mu}(t,Y)$ associated to an insider portfolio $\pi(t,Y)$ (representing the fraction of the wealth invested in the risky asset) is given by:
\begin{equation}\label{eq00.2}
\begin{cases}
dX(t,Y)=\pi(t,Y) X(t,Y) [ (\alpha(t,Y) +\mu(t))dt + \beta(t,Y)]dB(t) \\
X(0)=x >0
\end{cases}
\end{equation}
Define the performance functional by
\begin{equation}\label{eq00.3}
J(\pi,\mu)= \mathbb{E}[\int_0^T{\frac{1}{2} \mu^2(t)}dt + \theta(T,Y)\log X(T)],
\end{equation}
where $\theta(T,Y)>0$ is a given $\mathcal{H}_T$-measurable random variable, and $\frac{1}{2} \mu^2(t)$ represents a penalty rate, penalizing $\mu$ for being away from 0.
We assume that $\pi$ is $\mathbb{H}$-adapted, while $\mu$ is $\mathbb{F}$-adapted, i.e. has no inside information.

We want to find $\pi^* \in \mathcal{A}_1$ and $\mu^* \in \mathcal{A}_2$ such that

\begin{equation}\label{eq00.4}
\sup_{\pi \in \mathcal{A}_1} \inf_{\mu \in \mathcal{A}_2} J(\pi,\mu)=\inf_{\mu \in \mathcal{A}_2} \sup_{\pi \in \mathcal{A}_1} J(\pi,\mu)=J(\pi^*,\mu^*).
\end{equation}

We rewrite this problem as a classical stochastic differential game with one parameter $y_1=y \in \mathbb{R}$.
Thus we define
\begin{align}\label{eq00.5}
\begin{cases}
dx(t,y)&= \pi(t,y)x(t,y)[\{\alpha(t,y) +\mu(t)\}dt + \beta(t,y)dB(t)] \\
x(0,y)&=x(y) >0
\end{cases}
\end{align}
and
\begin{align}\label{eq00.6}
J(\pi(.,y),\mu(.))=\mathbb{E}[\int_0^T\ \frac{1}{2} \mu^2(t) \mathbb{E}[\delta_Y(y)|\mathcal{F}_t] dt + \theta \log x(T,y)\mathbb{E}[\delta_Y(y)|\mathcal{F}_T]].
\end{align}

The Hamiltonian for this problem is
\begin{align}\label{eq00.7}
H(t,x,y,\pi,\mu,p,q)= \frac{1}{2}\mu^2\mathbb{E}[\delta_Y(y)|\mathcal{F}_t] +\pi x (\alpha(t,y)+\mu)p +\pi x \beta(t,y)q
\end{align}
and the BSDE for the adjoint processes $p,q$ is

\begin{align}\label{eq00.8}
\begin{cases}
dp(t,y)= -[\pi(t,y)\{ (\alpha(t,y)+\mu(t))p(t,y)+\beta(t,y)q(t,y)\} ]dt
+ q(t,y)dB(t)  0 \leq t \leq T\\
p(T,y)=\frac{\theta(T,y) \mathbb{E}[\delta_Y(y)|\mathcal{F}_T]}{x(T,y)}.
\end{cases}
\end{align}

Maximizing $H$ with respect to $\pi$ gives the first order equation
\begin{equation}
x(t,y([(\alpha(t,y)+\mu(t))p(t,y) + \beta(t,y)q(t,y)]=0.
\end{equation}
Since $x(t,y) > 0$ and $\beta(t,y) \neq 0$, we deduce that
\begin{equation}
(\alpha(t,y)+\mu(t))p(t,y) + \beta(t,y)q(t,y)=0
\end{equation}
and
\begin{equation}
q(t,y)=-\frac{\alpha(t,y)+\mu(t)}{\beta(t,y)}p(t,y).
\end{equation}
Hence \eqref{eq00.6} reduces to
\begin{align}\label{eq00.12}
\begin{cases}
dp(t,y)= -\frac{\alpha(t,y)+\mu(t)}{\beta(t,y)}p(t,y) dB(t)\\
p(T,y)=\frac{\theta(T,y) \mathbb{E}[\delta_Y(y)|\mathcal{F}_T]}{x(T,y)}.
\end{cases}
\end{align}

Define
\begin{equation}\label{eq00.13}
h(t,y)=p(t,y)x(t,y).
\end{equation}
Then by the It\^ o formula we get
\begin{align}\label{eq00.14}
\begin{cases}
dh(t,y)&=(\pi(t,y)\beta(t,y) - \frac{\alpha(t,y)+\mu(t)}{\beta(t,y)})h(t,y)dB(t)\\
h(T,y)&=p(T,y)x(T,y)=\theta E[\delta_Y(y) | \mathcal{F}_T].
\end{cases}
\end{align}
This BSDE has the solution
\begin{equation}\label{eq00.15}
h(t,y) =  \theta E[\delta_Y(y) |\mathcal{F}_t].
\end{equation}
Moreover, by the generalized Clark-Ocone formula we have
\begin{equation}\label{eq00.16}
(\pi(t,y)\beta(t,y) - \frac{\alpha(t,y)+\mu(t)}{\beta(t,y)})h(t,y)=D_t h(t) = \theta E[D_t \delta_Y(y) | \mathcal{F}_t],
\end{equation}
from which we get the following expression for our candidate $\hat{\pi}(t,y)$ for the optimal portfolio
\begin{equation}\label{eq00.17}
\hat{\pi}(t,y)=\frac{\alpha(t,y)+\hat{\mu}(t)}{\beta^2(t,y)}+\frac{E[D_t\delta_Y(y) |\mathcal{F}_t]}{\beta(t,y)E[\delta_Y(y) | \mathcal{F}_t]},
\end{equation}
where $\hat{\mu}(t)$ is the corresponding candidate for the optimal perturbation.
\vskip 0.3cm

Minimizing $\int_{\mathbb{R}}Hdy$ with respect to $\mu$ gives the following first order equation for the optimal $\hat{\mu}(t)$:
\begin{align}
\int_{\mathbb{R}} \{ \hat{ \mu}(t)\mathbb{E}[\delta_Y(y)|\mathcal{F}_t] +\hat{\pi}(t,y)\hat{x}(t,y)\hat{p}(t,y)\}dy=0,
\end{align}
i.e.,
\begin{align}\label{eq00.19}
\hat{\mu}(t)= -\frac{\int_{\mathbb{R}}\hat{\pi}(t,y)\hat{x}(t,y)\hat{p}(t,y)dy}{\int_{\mathbb{R}}\hat{x}(t,y)\hat{p}(t,y)dy}= - \int_{\mathbb{R}} \hat{\pi}(t,y)E[\delta_Y(y) |\mathcal{F}_t] dy.
\end{align}

We can now verify that ($\hat{\pi},\hat{\mu})$ satisfies all the conditions of the sufficient maximum principle, and hence we conclude the following:

\begin{theorem}[Optimal portfolio for an insider under model uncertainty]
The saddle point $(\pi^*(t,Y),\mu^*(t))$, where $\pi^*(t,Y)=\pi^*(t,y) | _{y=Y}$, of the stochastic differential game \eqref{eq00.2} is given by
the solution of the following coupled system of equations
\begin{align}\label{eq00.20}
\pi^*(t,y)= \frac{\alpha(t,y)+\mu^*(t)}{\beta^2(t,y)}+\frac{E[D_t\delta_Y(y) |\mathcal{F}_t]}{\beta(t,y)E[\delta_Y(y) | \mathcal{F}_t]},
\end{align}
and
\begin{align}\label{eq00.21}
\mu^*(t)=  - \int_{\mathbb{R}} \pi^*(t,y)E[\delta_Y(y) |\mathcal{F}_t] dy.
\end{align}

\end{theorem}

\begin{remark}
This result is an extension to insider trading of a result in \cite{OS4}.
\end{remark}

 Consider the special case when $Y$ is a Gaussian random variable of the form
\begin{equation}\label{eqY}
    Y = Y (T_0); \text{ where } Y (t) =\int_0^t\psi(s)dB(s), \mbox{ for } t\in [0,T_0]
\end{equation}
for some deterministic function $\beta\in \mathbf{L}^2[0,T_0]$ with
\begin{equation}
    \|\psi\|^2_{[0,T]} :=\int_t^T\psi(s)^2ds>0 \mbox{ for all } t\in[0,T].
\end{equation}
In this case it is well known that the Donsker delta functional is given
by
\begin{equation}
    \delta_{Y}(y)=(2\pi v)^{-\frac{1}{2}}\exp^{\diamond}[-\frac{(Y-y)^{\diamond2}}{2v}]
\end{equation}
where we have put $v :=\|\psi\|^2_{[0,T_0]}$. See e.g. \cite{AaOU}, Proposition $3.2$.
We have 
\begin{equation}
\mathbb{E}[\delta_Y(y)|\mathcal{F}_t]=(2\pi \|\psi\|^2_{[t,T_0]})^{-\frac{1}{2}} \exp[- \frac{(Y(t)-y)^2}{2\|\psi\|^2_{[t,T_0]}}].
\end{equation}
and
\begin{equation}
\mathbb{E}[D_t\delta_Y(y)|\mathcal{F}_t] = -(2\pi \|\psi\|^2_{[t,T_0]})^{-\frac{1}{2}})\exp[- \frac{(Y(t)-y)^2}{2\|\psi\|^2_{[t,T_0]})}]\frac{Y(t)-y}{\|\psi\|^2_{[t,T_0]}}\psi(t).
\end{equation}
For more details see \cite{DrO}.
\begin{coro}
Suppose that $Y$ is Gaussian of the form (\ref{eqY}). Then the saddle point $(\pi^*(t,Y),\mu^*(t))$, where $\pi^*(t,Y)=\pi^*(t,y) | _{y=Y}$, of the stochastic differential game \eqref{eq00.2} is given by
the solution of the following coupled system of equations
\begin{align}\label{eq00.20}
\pi^*(t,Y)= \frac{\alpha(t,y)+\mu^*(t)}{\beta^2(t,y)}+\frac{Y(T_0)-Y(t)}{\beta(t,y) \|\psi\|^2_{[t,T_0]}}\psi(t),
\end{align}
and
\begin{align}\label{eq00.21}
\mu^*(t)=  -\theta (2\pi \|\psi\|^2_{[t,T_0]})^{-\frac{1}{2}}\int_{\mathbb{R}}\exp[- \frac{(Y(t)-y)^2}{2\|\psi\|^2_{[t,T_0]}}] \pi^*(t,y) dy.
\end{align}

\end{coro}
\begin{coro}
Suppose that $Y = B(T_0)$ for some $T_0 > T$. Then the saddle point $(\pi^*(t,Y),\mu^*(t))$, where $\pi^*(t,Y)=\pi^*(t,y) | _{y=Y}$, of the stochastic differential game \eqref{eq00.2} is given by
the solution of the following coupled system of equations
\begin{align}\label{eq00.20}
\pi^*(t,y)= \frac{\alpha(t,y)+\mu^*(t)}{\beta^2(t,y)}+\frac{B(T_0)-B(t)}{\beta(t,y) (T_0-t)}; \quad 0≤ t ≤T
\end{align}
and
\begin{align}\label{eq00.21}
\mu^*(t)=  -\theta (2\pi (T_0-t))^{-\frac{1}{2}}\int_{\mathbb{R}}\exp[- \frac{(B(t)-y)^2}{2(T_0-t)}] \pi^*(t,y) dy; \quad 0≤ t ≤T.
\end{align}
\end{coro}


\begin{thebibliography}{9999}

\bibitem[A\O ]{AO}
N. Agram and B. \O ksendal:
Malliavin calculus and optimal control of stochastic Volterra equations.
J. Optim. Theory Appl. (2015) DO1 10. 1007/s 10957-015-0753-5.

\bibitem[Aa\O PU]{AaOPU}
K. Aase, B. \O{}ksendal, N. Privault and J. Ub\o e: White noise generalizations of the Clark-Haussmann-Ocone theorem with application to mathematical finance. Finance Stoch. 4 (2000), 465-496.

\bibitem[Aa\O U]{AaOU}
K. Aase, B. \O{}ksendal and J. Ub\o e: Using the Donsker delta function to compute hedging strategies. Potential Analysis 14 (2001), 351-374.

\bibitem[B]{B} L. Breiman: Probability. Addison-Wesley 1968.

\bibitem[BBS]{BBS}
O. E. Barndorff-Nielsen, F.E. Benth and B. Szozda: On stochastic integration for volatility modulated Brownian-driven Volterra processes via white noise analysis. arXiv:1303.4625v1, 19 March 2013.

\bibitem[B\O]{BO}
F. Biagini and B. \O{}ksendal: A general stochastic calculus approach to insider trading.Appl. Math. \& Optim. 52 (2005), 167-181.

\bibitem [DM\O R]{DMOR}K.R. Dahl, S.-E. A. Mohammed, B. \O ksendal and E. R. R\o se; Optimal control with noisy memory and BSDEs with Malliavin derivatives. arXiv: 1403.4034 (2014).

\bibitem[DM\O P1]{DMOP1} G. Di Nunno, T. Meyer-Brandis, B. {\O}ksendal and F. Proske: Malliavin calculus and anticipative It\^ o formulae for L\'{e}vy processes. Inf. Dim. Anal. Quantum Prob. Rel. Topics 8 (2005), 235-258.

\bibitem[DM\O P2]{DMOP2} G. Di Nunno, T. Meyer-Brandis, B. {\O}ksendal and F. Proske: Optimal portfolio for an insider in a market driven by L\'{e}vy processes. Quant. Finance 6 (2006), 83-94.

\bibitem[D\O ]{DiO} G. Di Nunno and B. {\O}ksendal: The Donsker delta function, a representation formula for functionals of a L\'{e}vy process and application to hedging in incomplete markets. S\'{e}minaires et Congr\`{e}es, Societ\'{e} Math\'{e}matique de France, Vol. 16 (2007), 71-82.

\bibitem[D\O ]{DO} G. Di Nunno and B. \O ksendal: A representation theorem and a sensitivity result for functionals of jump diffusions. In A.B. Cruzeiro, H. Ouerdiane and N. Obata (editors): Mathematical Analysis and Random Phenomena. World Scientific 2007, pp. 177 - 190.

\bibitem[D\O P]{DOP} G. Di Nunno, B. {\O}ksendal and F. Proske: Malliavin Calculus for L\'{e}vy Processes with Applications to Finance. Universitext, Springer 2009.
    
\bibitem[Dr\O ]{DrO}  O. Draouil, B. {\O}ksendal: A Donsker delta functional approach to optimal insider
control and applications to finance. Communications in Mathematics and Statistics
(CIMS) 2015 (to appear). http://arxiv.org/abs/1504.02581.    

\bibitem[H\O UZ]{HOUZ} H. Holden, B. {\O}ksendal, J. Ub\o e and T. Zhang: Stochastic Partial Differential Equations. Universitext, Springer, Second Edition 2010.
    
\bibitem[LP]{LP} A. Lanconelli and F. Proske: On explicit strong solution of It\^ o-SDEs and the Donsker delta function of a diffusion.
Inf. Dim. Anal. Quatum Prob Rel. Topics 7 (2004),437-447.

\bibitem[M\O P]{MOP} S. Mataramvura, B. \O ksendal and F. Proske: The Donsker delta function of a L\'{e}vy process with application to chaos expansion of local time. Ann. Inst H. Poincar\'{e} Prob. Statist. 40 (2004), 553-567.

\bibitem[MP]{MP} T. Meyer-Brandis and F. Proske: On the existence and explicit representability of strong solutions of L\'{e}vy noise driven SDEs with irregular coefficients. Commun. Math. Sci. 4 (2006), 129-154.

\bibitem[\O R1]{OR1} B. \O ksendal and E. R\o se: A white noise approach to insider trading. Manuscript University of Oslo, 20 November 2014

\bibitem[\O R2]{OR2} B. \O ksendal and E. R\o se: Applications of white noise to mathematical finance. Manuscript University of Oslo, 5 February 2015

\bibitem[\O S1]{OS1}
B. \O{}ksendal and A. Sulem: Applied Stochastic Control of Jump Diffusions. Second Edition. Springer 2007

\bibitem[\O S2]{OS2}
B. \O{}ksendal and A. Sulem: Risk minimization in financial markets modeled by It\^ o-L\' evy processes. Afrika Matematika (2014), DOI: 10.1007/s13370-014-02489-9.

\bibitem[\O S3]{OS3}
B. \O{}ksendal and A. Sulem: A game theoretic approach to martingale measures in incomplete markets. Surveys of Applied and Industrial Mathematics, TVP Publishers, Moscow, 15 (2008), 18-24.

\bibitem[\O S4]{OS4}
B. \O{}ksendal and A. Sulem: Dynamic robust duality in utility maximization. ArXiv 1304.5040 (2014)

\bibitem[P]{P} P. Protter: Stochastic Integration and Differential Equations. Second Edition. Springer 2005

\bibitem[PK]{PK} I. Pikovsky and I. Karatzas: Anticipative portfolio optimization. Adv. Appl. Probab. 28 (1996), 1095-1122.

\bibitem[RV]{RV} F. Russo and P. Vallois: Forward, backward and symmetric stochastic integration. Probab. Theor. Rel. Fields 93 (1993), 403-421.

\bibitem[RV1]{RV1} F. Russo and P. Vallois. The generalized covariation process and It\^{o} formula. Stoch. Proc. Appl., 59(4):81-104, 1995.

\bibitem[RV2]{RV2} F. Russo and P. Vallois. Stochastic calculus with respect to continuous finite quadratic variation processes. Stoch. Stoch. Rep., 70(4):1-40, 2000.


\end{thebibliography}
\end{document}